\documentclass[a4paper, 11pt, babel]{amsart}

\usepackage{amsmath,enumerate, amsfonts, amssymb, amsthm, graphicx}
\usepackage[all]{xy}
\newtheorem{definition}{DEFINITION}[section]
\newtheorem{theorem}[definition]{THEOREM}
\newtheorem{proposition}[definition]{PROPOSITION}
\newtheorem{remark}[definition]{REMARK}
\newtheorem{lemma}[definition]{LEMMA}
\newtheorem{corollary}[definition]{COROLLARY}
\newtheorem{example}[definition]{EXAMPLE}

\input xy
\xyoption {all}
\newcommand{\sfrac}[2]{{#1}/{#2}}
\newcommand{\alg}{\mathrm{alg}}
\newcommand{\an}{\mathrm{an}}
 
\DeclareMathOperator{\spe}{sp} \DeclareMathOperator{\DR}{DR}

\title[Formal structure of direct image of holonomic $\mathcal{D}$-modules]{Formal structure of direct image of holonomic $\mathcal{D}$-modules of exponential type}
\author{C\'eline Roucairol}
\address{Lehrstuhl VI f\"ur Math., Univ. Mannheim,
A5 6, 68161 Mann\-heim, Germany.}
\email{celine.roucairol@uni-mannheim.de} 

\begin{document}

\begin{abstract}
We compute formal invariants associated with the cohomology
sheaves of the direct image of holonomic $\mathcal{D}$-modules of exponential
type. We also prove that every formal $\mathbb{C}[[t]]\langle\partial_t\rangle$-modules is isomorphic, after a ramification, to a germ of formalized direct image of analytic $\mathcal{D}$-module of exponential type.

\end{abstract}\maketitle
\section*{Introduction}
Let $X$ be a complex manifold. We denote by $\mathcal{O}_X$ (resp. $\mathcal{D}_X$) the sheaf of holomorphic functions (resp. holomorphic differential operators) on $X$. For any reduced divisor $Z$ in $X$, we denote by $\mathcal{O}_X[\ast Z]$ the sheaf of meromorphic functions with poles along $Z$ at most. Given a regular holonomic $\mathcal{D}_X$-module $\mathcal{M}$ and a meromorphic function $g$ on $X$ with poles along $Z$, we define the $\mathcal{D}_X$-module $\mathcal{M}e^g$ as being the $\mathcal{O}_X$-module $\mathcal{M}[\ast Z]:=\mathcal{O}_X[\ast Z]\otimes_{\mathcal{O}_X}\mathcal{M}$ with connection $\nabla_g$ defined as $\nabla+dg$, where $\nabla$ is the connection on $\mathcal{M}[\ast Z]$ given by its left $\mathcal{D}_X$-module structure. It is known that $\mathcal{M}e^g$ is a holonomic $\mathcal{D}_X$-module. It is irregular along $Z$. We will say $\mathcal{M}e^g$ has \emph{exponential type}.

In the algebraic setting, if $U$ is a smooth algebraic variety over $\mathbb{C}$, we use an exponent ``$\alg$'' in the notation of the corresponding sheaves and we define the notion of a holonomic $\mathcal{D}_X^{\alg}$-module of exponential type similarly.

Let $\mathcal{M}e^g$ be a holonomic $\mathcal{D}_X$-module of exponential type. The cohomology sheaves of its direct image by any proper map $f:X\to C$ to a complex curve are holonomic $\mathcal{D}_C$-modules with possibly irregular singularities. The main result of this article consists in the computation of the formal invariants of these cohomology modules at their singularities.

The basic computation takes place in the following setting: the complex manifold $X$ is the product $D\times \mathbb{P}^1$ of a small disc $D$ centered at the origin in $\mathbb{C}$ by the Riemann sphere, the meromorphic function $g$ is the projection $p_2:X\to\mathbb{P}^1$ with polar divisor $Z=D\times\{\infty\}$ and the map $f$ is the projection $p_1:X\to D$. We will denote by $t$ a coordinate on $D$. As above, $\mathcal{M}$ is a regular holonomic $\mathcal{D}_X$-module. The only interesting  cohomology module to consider is $\mathcal N=\mathcal{H}^0p_{1+}(\mathcal{M}e^{p_2})$ and we assume that $D$ is small enough so that $0$ is the only singular point of $\mathcal{N}$.

In a neighbourhood $U$ of $(0,\infty)\in D\times\mathbb{P}^1$, the singular support of $\mathcal{M}$ is an analytic curve $S$. We denote by $\{S_\ell\mid\ell\in\Lambda\}$ the set consisting of the local irreducible components of $S$ which are distinct of $\{0\}\times\mathbb{P}^1$ or $D\times\{\infty\}$. Set $U^{\ast}=U\setminus S$. To any $S_\ell$ at $(0,\infty)$ we associate:
\begin{figure}[h]
\includegraphics[height=2.4cm]{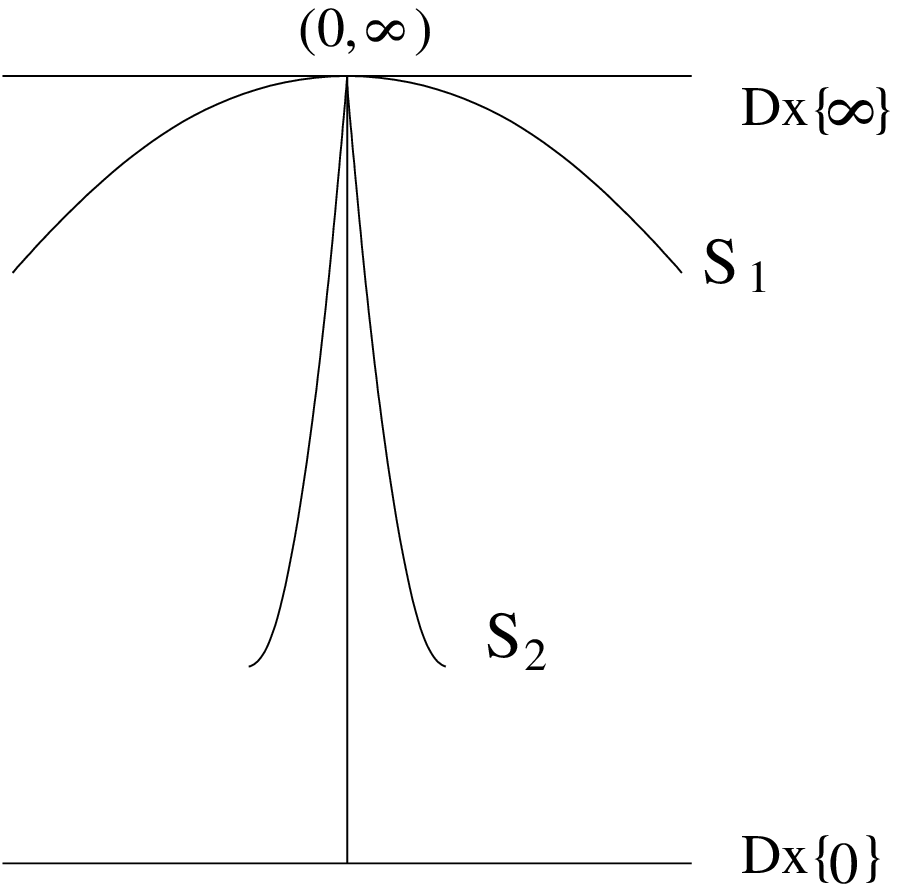}
\end{figure}
\begin{itemize}
\item
$m_\ell\in\mathbb{N}^*$:  the multiplicity of the conormal space $T^*_{S_\ell}X$ in the characteristic cycle of $\mathcal{M}$.
\item
$p_\ell,q_\ell\in\mathbb{N}^*$: the intersection multiplicity at $(0,\infty)$ of $S_\ell$ with $\{0\}\times \mathbb{P}^1$ and $D\times \{\infty\}$, respectively.
\item
$Q_\ell\subset\mathbb{R}^2$: the convex hull in $\mathbb{R}^2$ of the union of $Q:=\{(u,v)\in\mathbb{R}^2\mid u\leq0,\,v\geq0\}$ and $(m_\ell p_\ell,m_\ell q_\ell)+Q$.
\item $\alpha_\ell\in\tau^{-1}\mathbb{C}[\tau^{-1}]$ and $\delta_\ell\in\mathbb{C}\{\tau\}$: the polar part and the holomorphic part of a Puiseux parametrization of $S_\ell$ at $(0,\infty)$. Let $u_{\ell}\in\mathbb{C}\{\tau\}$, $u_\ell(0)\neq 0$, such that $\gamma(\tau):=(\tau^{p_\ell},\tau^{-q_\ell}u_\ell)$ is a parametrization of $S_\ell$. $\alpha_\ell$ is the polar part of $\tau^{-q_\ell}u_\ell$ and $\delta_\ell$ is its holomorphic part ($\tau^{-q_\ell}u_\ell=\alpha_\ell+\delta_\ell$).
\end{itemize}

The irregularity number of $\mathcal N_0$ is computed in [9]. But, we can go further in computing formal invariants of $\mathcal N_0$.

\begin{theorem}\label{1}\mbox{}\par
$(1)$ After a suitable translation, the Newton polygon of the $\mathcal D_{D,0}$-module $\mathcal N_0$ is the Minkowski sum $\sum_{\ell\in\Lambda}Q_\ell$. In particular, the set of slopes of $\mathcal N_0$ is $\{q_\ell/p_\ell\mid\ell\in\Lambda\}$ and the irregularity number of $\mathcal N_0$ is $\sum_{\ell\in\Lambda}m_\ell q_\ell$.

$(2)$ Let $p=\mathrm{lcm}\{p_\ell\mid\ell\in\Lambda\}$. After the ramification $\rho:D^{'}\to D$, $\rho(\tau)=t=\tau^p$, the formal irregular part of $\rho^*\mathcal N_0$ decomposes as $\oplus_{\alpha\in\Gamma}R_\alpha e^\alpha$ where:
\begin{itemize}
\item
$\Gamma\subset\tau^{-1}\mathbb C[\tau^{-1}]$ is a finite subset and $\alpha\in\Gamma$ if and only if there exist $\ell\in\Lambda$ and $\xi\in\mathbb{C}^{\ast}$ with $\xi^{p_\ell}=1$ such that
$$
\alpha(\tau)=\alpha_\ell(\xi\tau^{p/p_\ell});
$$
the set of such $\ell\in\Lambda$ is denoted by $\Lambda_\alpha$.
\item
$R_\alpha$ is a regular holonomic $\mathbb C[[\tau]]\langle\partial_\tau\rangle$-module, $R_\alpha=R_\alpha[\tau^{-1}]\neq0$ and the rank of $R_\alpha$ is $\sum_{\ell\in\Lambda_\alpha}m_\ell$.
\end{itemize}
\end{theorem}

But we can also compute the characteristic polynomial of the monodromy of $R_{\alpha}$ under the assumption $(\ast)$: for any $\ell,\ell^{'}\in\Lambda$ and for any $\xi_\ell,\xi_{\ell^{'}}\in\mathbb{C}^{\ast}$ with $\xi_\ell^{p_\ell}=1$ and $\xi_{\ell^{'}}^{p_{\ell^{'}}}=1$, 
$$(\ell,\xi_\ell)\neq(\ell^{'},\xi_{\ell^{'}})\Longrightarrow \alpha_\ell(\xi_\ell\tau^{\sfrac{p}{p_\ell}})+\delta_{\ell}(0)\neq\alpha_{\ell^{'}}(\xi_{\ell^{'}}\tau^{\sfrac{p}{p_{\ell^{'}}}})+\delta_{\ell^{'}}(0).$$

We remark that $\rho^{'}=(\rho,id):D^{'}\times\mathbb{P}^{1}\to D\times\mathbb{P}^1$ is a normalization of all the $S_{\ell}$'s. For any $\ell\in\Lambda$, $\rho^{'-1}(S_{\ell})$ is the union of some smooth irreductible analytic curves $S_{\ell}^i$.  

Then to any $S_{\ell}$ we associate $\zeta_{\ell}\in\mathbb{C}[\lambda]$: the characteristic polynomial of the monodromy of $\Phi_{h_{\ell}}(DR~\mathcal{M})$ around $(0,0)$ in a normalization of $S_{\ell}$, where $h_{\ell}=0$ is an equation of $S_{\ell}$. Here, $\Phi_{h_{\ell}}$ is the vanishing cycle functor along $S_{\ell}$ for bounded complex with constructible cohomology (cf. \cite{MeMai}, $\S 1$) and $\DR$ is the de Rham functor (cf. Definition $2.6.4$ p. $27$ in \cite{Me}). $\Phi_{h_{\ell}}(DR~\mathcal{M})$ has support included in $S_{\ell}$ and $\Phi_{h_{\ell}}(DR~\mathcal{M})_{|S_{\ell}\setminus\{(0,0)\}}$ is just a local system. Then, $\zeta_{\ell}$ is the characteristic polynomial of the monodromy of the local system  $\rho^{-1}(\Phi_{h_{\ell}}(DR~\mathcal{M}))_{|S_{\ell}^i\setminus\{(0,0)\}}$ around $(0,0)$ (do not confuse with the monodromy around $S_{\ell}$). As $S_{\ell}^i\simeq S_{\ell}^j$, for $i,j\in\{1,\ldots,p_{\ell}\}$, $\zeta_{\ell}$ is independant of $i$.

\begin{theorem}\label{zeta}
Under the assumption $(\ast)$, the characteristic polynomial of the monodromy of $R_\alpha$ is equal to $\prod_{\ell\in\Lambda_\alpha}\zeta_{\ell}$.
\end{theorem}

Conversely, we prove in Theorem \ref{realisation} that given any formal holonomic $\mathbb{C}[[t]]\langle\partial_t\rangle$-module $\widehat{\mathcal{N}}$, there exist a ramification $\rho:\tau\to t=\tau^p$ and a regular holonomic $\mathcal D_{D\times\mathbb P^1}$-module $\mathcal{M}$, such that $\rho^{\ast}(\widehat{\mathcal{N}})$ is isomorphic to the formalization of $\rho^{\ast}(\mathcal H^0p_{1+}(\mathcal{M}e^{p_2}))_0$.

In the algebraic setting, let us assume that $U$ is affine and that $f$ and $g$ are regular functions on $U$. Let $\mathcal{M}$ be a holonomic $\mathcal{D}^{\alg}_U$-module which has regular singularities included at infinity. The direct images $\mathcal{N}^k:=\mathcal{H}^k f_+(\mathcal{M}e^g)$ are holonomic $\mathbb C[t]\langle\partial_t\rangle$-modules with singularities at finite distance on the affine line $\mathbb{A}^1$ and at infinity. We will reduce the computation of the formal invariants of $\mathcal{N}^k$ at each of its singularities $c\in\mathbb{A}^1\cup\{\infty\}$ to the situation of Theorem \ref{1} through the diagram
\[
U\xrightarrow{~(f,g)~}\mathbb{A}^1\times\mathbb{A}^1\stackrel{i}{\to} \mathbb{P}^1\times\mathbb{P}^1\stackrel{i_c}{\leftarrow} D_c\times\mathbb{P}^1.
\]
Let $\mathcal{P}_c$ denote the direct image $i_c^{\ast}(i_+(f,g)_+\mathcal{M})^{\an}$. It consists in a complex of $\mathcal{D}_{D_c\times\mathbb{P}^1}$-modules with regular holonomic cohomology modules $\mathcal{H}^k\mathcal{P}_c$.
\begin{theorem}\label{2}
The germs $\mathcal{N}^k_c$ and $\mathcal{H}^0p_{1+}(\mathcal{H}^k\mathcal{P}_ce^{p_2})_c$ have the same formal irregular part.
\end{theorem}

\begin{figure}[h]
\includegraphics[height=2.5cm]{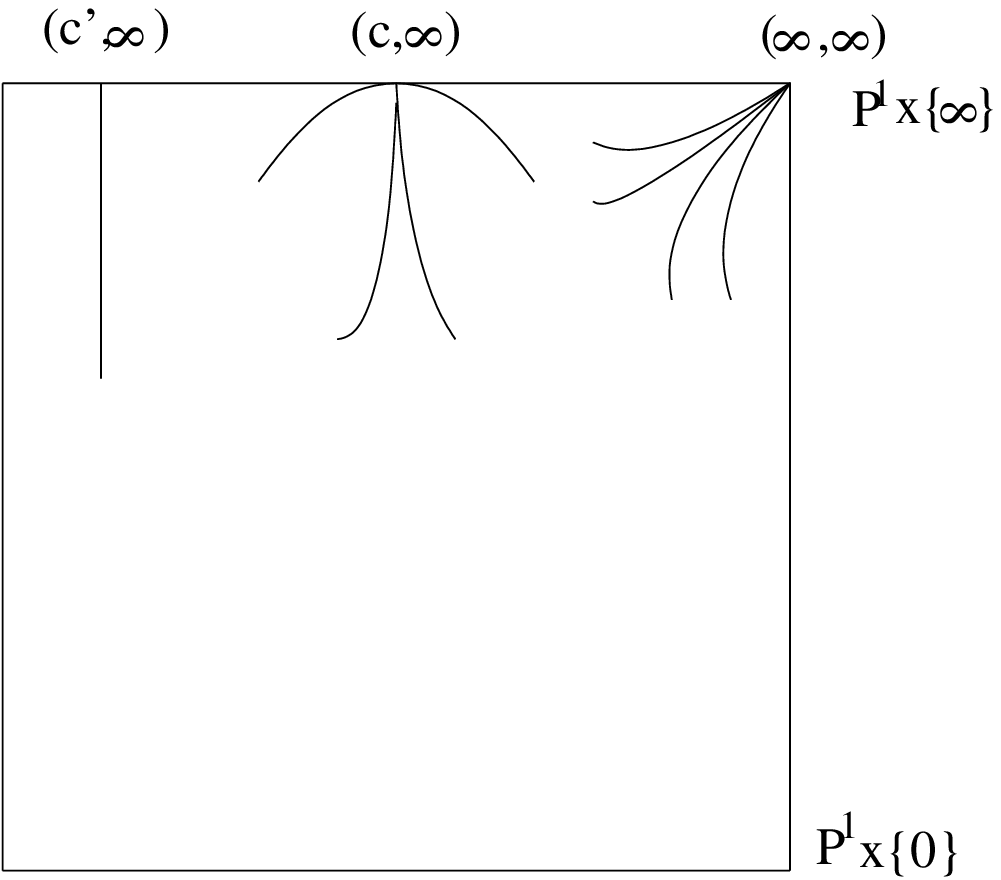}\\
Singular support of $\mathcal{H}^k(i_+(f,g)_+\mathcal{M})^{\an}$.
$c'$ is a regular singularity of $\mathcal{N}^k$. $c$ and $\infty$ are irregular singularities.
\end{figure}

\begin{example}
Let $U$ be a smooth affine surface. Suppose $f,g:U\to\mathbb{A}^1$ are algebraically independant and $\mathcal{M}=\mathcal{O}_U^{\alg}$. The only interesting cohomology module to consider is $\mathcal{N}=\mathcal{H}^0f_+(\mathcal{O}_U^{\alg}e^g)$; the others cohomology modules have punctual support (cf. Proposition 5.1 of \cite{Ro1}). 

Let $c\in\mathbb{A}^1\cup\{\infty\}$ and $S$ be the singular support of $\mathcal{H}^0i_+(f,g)_+\mathcal{O}_U$ in the neighbourhood of $(c,\infty)$. Theorem \ref{1} and Theorem \ref{2} allow us to describe the Newton polygon of $\mathcal{N}_c$, up to a suitable translation, using the set $\{S_{\ell}~|~\ell\in\Lambda\}$ of irreducible components of $S$  which are distinct of $\{c\}\times\mathbb{P}^1$ or $\mathbb{P}^1\times\{\infty\}$ and the multiplicity $m_{\ell}$ of $T^{\ast}_{S_{\ell}}(\mathbb{P}^1\times\mathbb{P}^1)$ in the characteristic cycle of $\mathcal{H}^0i_+(f,g)_+\mathcal{O}_U$.

But we can give a geometrical interpretation of the $S_{\ell}$'s and the $m_{\ell}$'s.
Let $\mathbb{X}$ be a smooth compactification of $U$ such that there exists two holomorphic functions $F,G:\mathbb{X}\to\mathbb{P}^1$  which extend $f$ and $g$. Let $\Gamma$ be the critical locus of $(F,G)$ and $D$ be the divisor $\mathbb{X}\setminus U$. 
We denote by $\Delta_1$ (resp. $\Delta_2$) the cycle in $\mathbb{P}^1\times\mathbb{P}^1$ which is the closure of $(F,G)(\Gamma)\cap (\mathbb{A}^1\times\mathbb{A}^1)$ (resp. $(F,G)(D)\cap (\mathbb{A}^1\times\mathbb{A}^1)$), where the image is counted with multiplicity. If $c\in\mathbb{A}^1$ (resp. $c=\infty$), the $S_{\ell}$'s are the irreducible components of the germ $\Delta_1\setminus(\{c\}\times\mathbb{P}^1)$ (resp. $\Delta_1\cup\Delta_2$) at $(c,\infty)$; $m_{\ell}$ is the multiplicity of $S_{\ell}$ in $\Delta_1\setminus(\{c\}\times\mathbb{P}^1)$ (resp. $\Delta_1\cup\Delta_2$). This computation can be found in \cite{Ro2} (Theorem 6.4.1).
\end{example}

We shall assume that the reader is familiar with the definitions and the properties of the specialization functor and the nearby cycles functor for $\mathcal D$-modules along a hypersurface and we refer to \cite{Ka}, \cite{Ma2}, \cite{MeSa}, \cite{LaMa} and \cite{MeMai} for details. Given $Y\subset X$ a smooth hypersurface of $X$ defined by a global equation $h=0$ and given $\mathcal{M}$ a specializable $\mathcal{D}_X$-module, we will denote by $\spe_Y\mathcal{M}$ the specialization of $\mathcal{M}$ along $Y$ and by $\Psi_h\mathcal{M}$ the nearby cycle module of $\mathcal{M}$ along $Y$. If $V_{\bullet}\mathcal{M}$ is the canonical V-filtration with respect to the lexical order on $\mathbb{C}\simeq\mathbb{R}\oplus i\mathbb{R}$, $\spe_Y\mathcal{M}=\oplus_{\beta\in\mathbb{C}}\sfrac{V_{\beta}\mathcal{M}}{V_{<\beta}\mathcal{M}}$ and 
$\Psi_h\mathcal{M}=\oplus_{-1\leq\beta<0}\sfrac{V_{\beta}\mathcal{M}}{V_{<\beta}\mathcal{M}}$. Same definitions hold for $Y$  hypersurface with normal crossings.
We will also denote by $\Psi_h$ the nearby cycle functor for bounded complexes with constructible cohomology (cf \cite{MeMai} \S 1).

Let us explain the method we use to compute the formal invariants. We note that if a summand $R_\alpha e^\alpha$ shows up in the formal decomposition of a $\mathbb{C}\{t\}\langle\partial_t\rangle$-module $\mathcal{N}_0$ after a ramification $\rho$, then $R_\alpha$ is the formal regular part of $((\rho^*\mathcal{N})e^{-\alpha})_0$. This formal regular part can be recovered from the specialization of $(\rho^*\mathcal{N})e^{-\alpha}$ at the origin (cf. Example 5.2.1 in \cite{LaMa}). Moreover, as $R_{\alpha}$ is localized at the origin, the rank of $R_{\alpha}$ and the characteristic polynomial of the monodromy of $R_{\alpha}$ can be computed using the nearby cycles module of $(\rho^*\mathcal{N})e^{-\alpha}$ at the origin.
One of the main results we use is the commutation of the specialization functor (resp. nearby cycles functor) with proper direct images (cf. Theorem $9.4.1$ in \cite{LaMa} for complexes of $\mathcal{D}$-modules and Theorem $4.8-1$ in \cite{MeSa} for $\mathcal{D}$-modules). This property reduces the proof of Theorem \ref{1} to local analytic computations on a suitable blowing-up space of $D\times\mathbb P^1$ above $(0,\infty)$. It is also the main argument in the proof of Theorem \ref{2}.

\section{Direct image of holonomic $\mathcal{D}$-module of exponential type: the case of projections}

This section is devoted to the proof of Theorems \ref{1} and \ref{zeta}. We first note that Theorem \ref{1} $(2)$ induces the construction of the Newton polygon (Theorem \ref{1} $(1)$).

Let $s\in\mathbb{N}^{\ast}$. Denote by $\Gamma_s$ the set of $\alpha\in\Gamma$ with pole of order $s$ and by $\Lambda_s$ the set of $\ell\in\Lambda$ such that $\sfrac{pq_{\ell}}{p_{\ell}}=s$. The Newton polygon of $\oplus_{\alpha\in\Gamma_s}R_{\alpha}e^{\alpha}$ has just one slope, equal to $s$; its height is the product of the rank of $\oplus_{\alpha\in\Gamma_s}R_{\alpha}$ by the slope. According to Theorem \ref{1} $(2)$, the slope is equal to $\sfrac{pq_{\ell}}{p_{\ell}}$, for all $\ell\in\Lambda_s$, and the height is  $\sum_{\alpha\in\Gamma_s}\sum_{\ell\in\Lambda_{\alpha}}\sfrac{pq_{\ell}m_{\ell}}{p_{\ell}}=\sum_{\ell\in\Lambda_s}pq_{\ell}m_{\ell}$. If $\widetilde{Q_{\ell}}$ denotes the convex hull in $\mathbb{R}^2$ of the union of $Q$ and $(m_{\ell}p_{\ell},pm_{\ell}q_{\ell})+Q$, the Newton polygon of $\rho^{\ast}\mathcal{N}_0$ is, after a suitable translation, the Minkowski sum $\sum_{\ell\in\Lambda}\widetilde{Q_{\ell}}$. We deduce the Newton polygon of $\mathcal{N}_0$ from the one of $\rho^{\ast}\mathcal{N}_0$ by a dilation of the vertical axis in a ratio $\sfrac{1}{p}$ (cf. Lemma 5.4.3 p. $34$ of \cite{Sa}).

We focus now on the proof of Theorem \ref{1} $(2)$ and Theorem \ref{zeta}.

\subsection{Ramification}

The choice of the ramification enables us to reduce the proof of Theorem \ref{1} and \ref{zeta} to the case where the $p_l$'s are equal to $1$; it is the non ramified case ($p=1$).

We begin with a base change formula. Let $p_1^{'}:D^{'}\times\mathbb{P}^1\to D^{'}$ and $p_2^{'}:D^{'}\times\mathbb{P}^1\to\mathbb{P}^1$ be the canonical projections. Consider the cartesian diagram:
$$\xymatrix{
D^{'}\times\mathbb{P}^1\ar[r]^{\rho^{'}}\ar[d]^{p_1^{'}}&D\times\mathbb{P}^1\ar[d]^{p_1}\\
D^{'}\ar[r]^{\rho}&D.}$$

\begin{lemma}\label{comm}
$\rho^{\ast}p_{1+}(\mathcal{M}e^{p_2})=p_{1+}^{'}(\rho^{'\ast}(\mathcal{M})e^{p_2^{'}})$.
\end{lemma}

\begin{proof}
Denote by $~^p\DR_{D\times\mathbb{P}^1/D}$ (resp. $~^p\DR_{D^{'}\times\mathbb{P}^1/D^{'}}$) the relative de Rham functor of $p_1$ (resp. $p_1^{'}$). We adopt the convention that the relative de Rham complexes are concentrated in negative degrees. We have isomorphisms of complexes of $\mathcal{D}_{D^{'}}$-modules:
$$
\begin{array}{ll}
\rho^{\ast}p_{1+}&(\mathcal{M}e^{p_2})=\\
&=\mathcal{D}_{D^{'}\to D}\otimes_{\rho^{-1}\mathcal{D}_D}\rho^{-1}\mathbb{R}p_{1\ast}
~^p\DR_{D\times\mathbb{P}^1/D}(\mathcal{M}e^{p_2})\text{ (by definition),}\\
&=\mathcal{D}_{D^{'}\to D}\otimes_{\rho^{-1}\mathcal{D}_D}\mathbb{R}p^{'}_{1\ast}\rho^{'-1}
~^p\DR_{D\times\mathbb{P}^1/D}(\mathcal{M}e^{p_2}),\\
&\text{(Proposition 2.6.7 in \cite{KaSc}),}\\
&=\mathbb{R}p^{'}_{1\ast}(\mathcal{D}_{D^{'}\times\mathbb{P}^1\to D\times\mathbb{P}^1}\otimes_{\rho^{'-1}\mathcal{D}_{D\times\mathbb{P}^1}}\rho^{'-1}
~^p\DR_{D\times\mathbb{P}^1/D}(\mathcal{M}e^{p_2})),\\
&\text{ (Proposition 2.6.6 in \cite{KaSc}),}\\
&=\mathbb{R}p^{'}_{1\ast}~^p\DR_{D^{'}\times\mathbb{P}^1/D^{'}}((\mathcal{D}_{D^{'}\times\mathbb{P}^1\to D\times\mathbb{P}^1}\otimes_{\rho^{'-1}\mathcal{D}_{D\times\mathbb{P}^1}}\rho^{'-1}
(\mathcal{M}))e^{p_2^{'}}),\\
&=p_{1+}^{'}(\rho^{'\ast}(\mathcal{M})e^{p_2^{'}}).
\end{array}$$

\end{proof}

We deduce that $\rho^{\ast}\mathcal{N}=\mathcal{H}^0p_{1+}^{'}(\rho^{'\ast}(\mathcal{M})e^{p_2^{'}})$.

Then we remark that the set $\{S_{\ell}^i~|~\ell\in\Lambda,~i=1,\ldots,p_{\ell}\}$ consists of the irreducible components of the singular support $S^{'}$ of $\rho^{'\ast}(\mathcal{M})$ in a neighbourhood $U^{'}$ of $(0,\infty)$, which are distinct of $\{0\}\times\mathbb{P}^1$ or $D^{'}\times\{\infty\}$. Let $\xi_1,\ldots,\xi_{p_{\ell}}$ be the $p_{\ell}$-roots of the unity. Up to change of indices $i$, a parametrization of $S_{\ell}^i$ is given by $\gamma(\tau)=(\tau,\alpha_{\ell}(\xi_i\tau^{\sfrac{p}{p_{\ell}}})+\delta_{\ell}(\xi_i\tau^{\sfrac{p}{p_{\ell}}}))$; in particular, the intersection multiplicity 
$p_{\ell}^i$ at $(0,\infty)$ of $S_{\ell}^i$ with $\{0\}\times\mathbb{P}^1$ is $1$.

In the next section, we will prove Theorem \ref{1} $(2)$ and Theorem \ref{zeta} in the case where all the $p_{\ell}$'s are equal to $1$. Then, we may apply it to $\rho^{'\ast}(\mathcal{M})$. At last, we have to compute the data associated with each $S_{\ell}^i$ using data associated with the $S_{\ell}$'s.
But we have:\\

$\bullet$ the polar part of a Puiseux parametrization of $S_{\ell}^i$ at $(0,\infty)$ is $\alpha_{\ell}(\xi_i\tau^{\sfrac{p}{p_{\ell}}})$. Moreover, if the assumption $(\ast)$ is fulfilled for $\mathcal{M}$ it is also fulfilled for $\rho^{'\ast}(\mathcal{M})$, \\

$\bullet$ $m_{\ell}$ is the multiplicity of  $T^{\ast}_{S_{\ell}^i}(D^{'}\times\mathbb{P}^1)$ in the characterstic cycle of $\rho^{'\ast}(\mathcal{M})$,\\

$\bullet$ $\zeta_{\ell}$ is the characteristic polynomial of the monodromy of the local system  $\Phi_{h_{\ell}^i}\DR~(\rho^{'\ast}(\mathcal{M}))_{|S_{\ell}^i\setminus\{(0,0)\}}$ around $(0,0)$, where $h_{\ell}^i=0$ is an equation of $S_{\ell}^i$.

\subsection{Non ramified case: $p_{\ell}=1$, for all $\ell\in\Lambda$}
If we prove that the rank of $R_{\alpha}$ is equal to $\sum_{\ell\in\Lambda_{\alpha}}m_{\ell}$, the irregularity number at $0$ of $\mathcal{N}$ and of $\oplus_{\alpha\in\Gamma}R_{\alpha}e^{\alpha}$ are both equal to $\sum_{\ell\in\Lambda}m_{\ell}q_{\ell}$ (cf. \cite{Ro} for $\mathcal{N}$). Then it is sufficient to prove that each summand $R_{\alpha}e^{\alpha}$, with $\alpha\in\Gamma$, shows up in the formal decomposition of $\mathcal{N}_0$ and to compute the rank and the characteristic polynomial of the monodromy of each $R_{\alpha}$.

Let $\alpha\in\Gamma$. We denote by $\spe_0(\mathcal{N}e^{-\alpha})$ the specialization of $\mathcal{N}e^{-\alpha}$ at the origin. $\spe_0(\mathcal{N}e^{-\alpha})_0$ is a $\mathbb{C}[t]\langle\partial_t\rangle$-module. As announced in the introduction, $R_{\alpha}$ is the formal regular part of $(\mathcal{N}e^{-\alpha})_0$. It can be recovered using the isomorphism $R_{\alpha}=\mathbb{C}[[t]]\otimes_{\mathbb{C}[t]}\spe_0(\mathcal{N}e^{-\alpha})_0$ (cf. Example $5.2.1$ of \cite{LaMa}).
Let $\Psi_t(\mathcal{N}e^{-\alpha})$ be the nearby cycles module of $\mathcal{N}e^{-\alpha}$ at the origin. $\Psi_t(\mathcal{N}e^{-\alpha})_0$ is a finite dimensional $\mathbb{C}$-vector space equipped with an endomorphism of monodromy induced by $\exp(-2i\pi t\partial_t)$.
As $\mathcal{N}e^{-\alpha}$ is localized at the origin, the rank of $R_{\alpha}$ is the dimension of $\Psi_t(\mathcal{N}e^{-\alpha})_0$ and the characteristic polynomials of their monodromies are equal.

But we have, $\Psi_t(\mathcal{N}e^{-\alpha})_0=R\Gamma(\{0\}\times\mathbb{P}^1,\DR~\Psi_{p_1}(\mathcal{M}e^{p_2-\alpha\circ p_1})[+1])$. Indeed,

$$\begin{array}{ll}
\mathcal{N}e^{-\alpha}&=\mathcal{O}_D[\frac{1}{t}]e^{-\alpha}\otimes_{\mathcal{O}_D}\mathcal{H}^0p_{1+}(\mathcal{M}e^{p_2}),\\
&=\mathcal{O}_D[\frac{1}{t}]e^{-\alpha}\otimes^{\mathbb{L}}_{\mathcal{O}_D}\mathcal{H}^0p_{1+}(\mathcal{M}e^{p_2}),\text{ (flatness),}\\
&=(\mathcal{O}_D[\frac{1}{t}]e^{-\alpha}\otimes_{\mathcal{O}_D}\mathcal{D}_D)\otimes^{\mathbb{L}}_{\mathcal{D}_D}\mathcal{H}^0p_{1+}(\mathcal{M}e^{p_2}),\\
&=\mathcal{H}^0p_{1+}((\mathcal{O}_{D\times\mathbb{P}^1}[\frac{1}{t}]e^{-\alpha\circ p_1}\otimes_{\mathcal{O}_{D\times\mathbb{P}^1}}\mathcal{D}_{D\times\mathbb{P}^1})\otimes_{\mathcal{D}_D}^{\mathbb{L}}\mathcal{M}e^{p_2}),\\
&\text{ (Prop 2.6.6 in \cite{KaSc}),}\\
&=\mathcal{H}^0p_{1+}(\mathcal{M}e^{p_2-\alpha\circ p_1}).
\end{array}$$

Let $\overline{p_1}:\{0\}\times\mathbb{P}^1\to\{0\}$. According to the commutation of the nearby cycle functor with proper direct image (cf. Theorem 4.8-1 p. $226$ of \cite{MeSa}),
$$\begin{array}{ll}
\Psi_t(\mathcal{N}e^{-\alpha})_0&=\mathcal{H}^0\overline{p_1}_+(\Psi_{p_1}(\mathcal{M}e^{p_2-\alpha\circ p_1}))_0,\\
&=R\Gamma(\{0\}\times\mathbb{P}^1,\DR~\Psi_{p_1}(\mathcal{M}e^{p_2-\alpha\circ p_1})[+1]), \text{ (cf. \cite{Ma3} p.$5$).}
\end{array}$$

Then, Theorem \ref{1} $(2)$ and Theorem \ref{zeta} follows from the proposition:
\begin{proposition}\label{11}
$\Psi_{p_1}(\mathcal{M}e^{p_2-\alpha\circ p_1})$ has support included in $(0,\infty)$ and the Euler characteristic of  $\DR~\Psi_{p_1}(\mathcal{M}e^{p_2-\alpha\circ p_1})_{(0,\infty)}$ is $-\sum_{\ell\in\Lambda_{\alpha}}m_{\ell}$.\\
Moreover, under the assumption $(\ast)$, the zeta function of its monodromy is $\prod_{\ell\in\Lambda_\alpha}\zeta_{\ell}^{-1}$.
\end{proposition}

The proof is given in the next sections.

\begin{remark}
In the non ramified case, the assumption $(\ast)$ is equivalent to:
$$\forall\ell,\ell^{'}\in\Lambda_{\alpha},~\ell\neq\ell^{'}\Longrightarrow\delta_\ell(0)\neq\delta_{\ell^{'}}(0).$$
\end{remark}

\subsection{Local computations of nearby cycles modules}
In this section, we prepare the proof of Proposition \ref{11} by giving some computations of nearby cycles modules along a normal crossing. Let $V$ be a neighbourhood of $(0,0)$ in $\mathbb{C}^2$ and $(u,v)$ be some coordinates on $V$. Let $\mathcal{M}$ be a regular holonomic $\mathcal{D}_V$-module and $m,n,k,l\in\mathbb{N}$. We are interested in nearby cycles modules of the type $\Psi_{u^mv^n}(\mathcal{M}e^{\sfrac{1}{u^kv^l}})$.

\subsubsection{Nearby cycles module along $u^mv^n=0$ of a module of exponential type}\label{nota}
Suppose that the singular support of $\mathcal{M}$ is included in $uv=0$. Let $V^{\ast}=V\setminus\{uv=0\}$. $\mathcal{M}_{|V^{\ast}}$ is a holomorphic connection; we denote by $r$ its rank. Thus $(\DR~\mathcal{M})_{|V^{\ast}}$ is a local system. For $P\in V^{\ast}$, we consider the monodromy $T_P$ of $(\DR~\mathcal{M})_P$ around $u=0$. As $V^{\ast}$ is connected and the singular support of $\mathcal{M}$ is included in a normal crossing, the characteristic polynomial of $T_P$ does not depend on $P$ (cf. I.2.2 p. 55 in \cite{MeNa}). We denote it by $\zeta_r$.

\begin{lemma}\label{B}
\begin{enumerate}
\item If $k,m\geq1$, $\Psi_{u^m}(\mathcal{M}e^{\sfrac{1}{u^k}})=0$.
\item If $k,l,m,\geq1$ and $n\geq 0$, $\Psi_{u^mv^n}(\mathcal{M}e^{\sfrac{1}{u^kv^l}})=0$.
\item If $n\geq 1$, the Euler characteristic of  $\DR~\Psi_{uv^n}(\mathcal{M}[\frac{1}{uv}]e^{\sfrac{1}{v}})_{(0,0)}$ is equal to $-r$ and the zeta function of its monodromy is $\zeta_r^{-1}$.
\end{enumerate}
\end{lemma}

\begin{proof}
The proof of the last two points is given in Lemma $4.5.10$ of \cite{Sa2} and its proof. The shift in the third point comes from the convention on the de Rham complexes. The first point is also proved in \cite{Sa2} in the case where $\mathcal{M}=\mathcal{M}[\frac{1}{v}]$. To conclude, it is sufficient to prove it in the case where $\mathcal{M}$ has support in $v=0$. We have $\mathcal{M}=i_+Li^{\ast}(\mathcal{M})[+1]$, where $i:\{v=0\}\hookrightarrow V$,
and $Li^{\ast}(\mathcal{M})[+1]$ is just a module.
As $\Psi_{u^m}(Li^{\ast}(\mathcal{M})[+1]e^{\sfrac{1}{u^k}}))=0$ (computation in one variable), we deduce that $\Psi_{u^m}(\mathcal{M}e^{\sfrac{1}{u^k}})=0$.
\end{proof}

\subsubsection{Nearby cycles module along $u=0$ of a regular $\mathcal{D}_V$-module}

Let $\widetilde{S}$ be the singular support of $\mathcal{M}$ and $\{\widetilde{S_{\ell}}~|~\ell\in\Lambda\}$ be the set consisting of the irreducible components of $\widetilde{S}$ which are distinct of $u=0$. We will assume that {\bf the intersection
multiplicity at $(0,0)$ of $\widetilde{S_{\ell}}$ with $u=0$ is equal to $1$}. 

We denote by $r$ the rank of $\mathcal{M}$. To any $\widetilde{S_{\ell}}$, we associate $m_{\ell}$, the multiplicity of the conormal space $T^{\ast}_{\widetilde{S_{\ell}}}V$ in the characteristic cycle of $\mathcal{M}$.

\begin{lemma}\label{A}
The Euler characteristic of $\DR~\Psi_u(\mathcal{M}[\frac{1}{u}])_{(0,0)}$ is equal to  $r-\sum_{\ell\in\Lambda}m_{\ell}$.
\end{lemma}

\begin{proof}
As $\mathcal{M}$ is regular holonomic, 
$$\begin{array}{ll}\DR~\Psi_u(\mathcal{M}[\frac{1}{u}])&=\DR~\Psi_u(\mathcal{M})\text{ (Proposition 2.4-3 in \cite{MeMai}),}\\
&=\Psi_u(\DR~\mathcal{M})\text{ (Theorem $4.10.1$ p. 233 in \cite{MeSa}).}\end{array}$$

Let $\epsilon$ and $\eta$ small enough such that any $\widetilde{S_{\ell}}$ intersects
$X_{\epsilon,\eta}=B(0,\epsilon)\cap\{u=\eta\}$ at a unique point $P_{\ell}$ (by assumption) and 
$\Psi_u(\DR~\mathcal{M})_{(0,0)}=\mathbb{R}\Gamma(X_{\epsilon,\eta},\DR~\mathcal{M})$ (cf. Remark $1.1-7$ in \cite{MeMai}).
Denote by
$X_{\epsilon,\eta}^c=X_{\epsilon,\eta}\setminus\{P_{\ell},~\ell\in\Lambda\}$.  Let $B_{\ell}$, $\ell\in\Lambda$, be some disjoint balls centered at $P_{\ell}$ of radius small enough. 

 According to the index theorem of Kashiwara (Theorem 6.3.1 in \cite{Ka1}), the Euler characteristic of $(\DR~\mathcal{M})_P$, $P\in X_{\epsilon,\eta}^c$, (resp. $(\DR~\mathcal{M})_{P_{\ell}}$) is equal to $r$ (resp. $r-m_{\ell}$).
As $(\DR~\mathcal{M})_{|X_{\epsilon,\eta}}$ is constructible
with respect to the stratification $\{X_{\epsilon,\eta}^c,
P_{\ell},~\ell\in\Lambda\}$, the Euler characteristic of  $\mathbb{R}\Gamma(X_{\epsilon,\eta},\DR~\mathcal{M})$ may be computed using the Mayer-Vietoris Theorem. We obtain:
$$\begin{array}{lll}
\chi(\mathbb{R}\Gamma(X_{\epsilon,\eta},\DR~\mathcal{M}))&=&
\chi(\mathbb{R}\Gamma(X_{\epsilon,\eta}^c,\DR~\mathcal{M}))+\sum_{\ell\in\Lambda}\chi((\DR~\mathcal{M})_{P_{\ell}})\\
&&-\sum_{\ell\in\Lambda}\chi(\mathbb{R}\Gamma(B_{\ell}\setminus\{P_{\ell}\},\DR~\mathcal{M})),\\
&=&\chi(X_{\epsilon,\eta}^c).r+\sum_{\ell\in\Lambda}(r-m_{\ell})\\
&&-\sum_{\ell\in\Lambda}\chi(B_{\ell}\setminus\{P_{\ell}\}).r,\\
&=&r-\sum_{\ell\in\Lambda}m_{\ell}.
\end{array}$$
\end{proof}

The zeta function of the monodromy of $\DR~\Psi_u(\mathcal{M}[\frac{1}{u}])_{(0,0)}$ can not be compute with this method in the general case. Nevertheless, we can do it when the singular support of $\mathcal{M}$ is included in the normal crossing $uv=0$. We adopt the notations of the section \ref{nota}. Let $V_0=\{v=0\}\setminus\{(0,0)\}$. We denote by $\zeta$ the characteristic polynomial of the monodromy of the local system $\phi_v(\DR~\mathcal{M})_{|V_0}$ around $(0,0)$.

\begin{lemma}\label{C}
If the singular support of $\mathcal{M}$ is included in a normal crossing $uv=0$, the zeta function of the monodromy of $\DR~\Psi_u(\mathcal{M}[\frac{1}{u}])_{(0,0)}$ is equal to $\zeta_r\zeta^{-1}$.
\end{lemma}

\begin{proof}
We have $\DR~\Psi_u(\mathcal{M}[\frac{1}{u}])_{(0,0)}=\underset{\underset{\epsilon,\eta}{\leftarrow}}{\lim}R\Gamma(X_{\epsilon,\eta},\DR~\mathcal{M})$.

As the singular support of $\mathcal{M}$ is included in $uv=0$ and as $\DR~\mathcal{M}$ is constructible with respect to a stratification such that $uv=0$ is a union of strata, we have $R\Gamma(X_{\epsilon,\eta},\DR~\mathcal{M})=(\DR~\mathcal{M})_P$, where $P\in V_0$ is the intersection point of $X_{\epsilon,\eta}$ with $v=0$. Then, the zeta function of the monodromy of $\DR~\Psi_u(\mathcal{M}[\frac{1}{u}])_{(0,0)}$ coincide with the zeta function of the monodromy of the complex of local systems $(\DR~\mathcal{M})_{|V_0}$ around $(0,0)$.

Consider the triangle of complexes of local systems:
$$(\DR~\mathcal{M})_{|V_0}\to\Psi_v(\DR~\mathcal{M})_{|V_0}\to\phi_v(\DR~\mathcal{M})_{|V_0}\overset{+1}{\to}.$$

As the zeta function of the monodromy of $\Psi_v(\DR~\mathcal{M})_{|V_0}$ around $(0,0)$ is $\zeta_r$, Lemma \ref{C} follows.
\end{proof}

\subsection{Proof of Proposition \ref{11}}~\\
First we have to prove that $\Psi_{p_1}(\mathcal{M}e^{p_2-\alpha\circ p_1})_{(0,c)}=0$, for any $c\neq \infty$. But in a neighbourhood $U$ of $(0,c)$ in $D\times\mathbb{P}^1$, there exist coordinates $(x,y)$ such that $p_1(x,y)=x$, $p_2(x,y)=y$ and $\alpha\circ p_1\sim\frac{1}{x^q}$, for some $q\in\mathbb{N}^{\ast}$. Therefore,  $\Psi_{p_1}(\mathcal{M}e^{p_2-\alpha\circ p_1})_{(0,c)}=\Psi_x(\mathcal{M}_{|U}e^{\frac{1}{x^q}})_{(0,0)}$. As the singular support of $\mathcal{M}_{|U}$ is not necessarily a normal crossing at $(0,0)$, we can not apply directly Lemma \ref{B} 1.

The idea is to use a resolution  $\pi=(\pi_1,\pi_2):\mathbb{X}\to U$, with exceptional locus $E$, such that the singular support of $\pi^{\ast}(\mathcal{M}_{|U})[\ast\pi_1^{-1}(0)]$ has normal crossings in a neighbourhood of any point of $E$. As $\mathcal{M}$ is holonomic and $\pi$ is a proper isomorphism out of $E$, 
$\mathcal{M}_{|U}e^{-\frac{1}{x^q}}=\pi_+(\pi^{\ast}(\mathcal{M}_{|U})e^{-\frac{1}{\pi_1^q}})$ (cf. Proposition 7.4.5 of \cite{Me}).
Using the commutation of the proper direct image with $\Psi$ and $\DR$ (Theorem 4.8.1 p. 226 in \cite{MeSa} and Theorem 5.4.3 p. 77 in \cite{Me}), we deduce that 
\label{9}
$$(\ast\ast)~\DR ~\Psi_{x}(\mathcal{M}_{|U}e^{-\frac{1}{x^q}})_{(0,0)}=R\Gamma(E,\DR~\Psi_{\pi_1}(\pi^{\ast}(\mathcal{M}_{|U})e^{-\frac{1}{\pi_1^q}})).$$
\begin{itemize}
\item If $P$ is a smooth point of $\pi^{-1}_1(0)$, we have $\Psi_{\pi_1}(\pi^{\ast}(\mathcal{M}_{|U})e^{-\frac{1}{\pi_1^q}})_P=0$ (cf. Lemma \ref{B} 1). 
\item If $P$ is a normal crossing point of $\pi^{-1}_1(0)$, according to Lemma \ref{B} 2, $\Psi_{\pi_1}(\pi^{\ast}(\mathcal{M}_{|U})e^{-\frac{1}{\pi_1^q}})_P=0$. 
\end{itemize}

Then $\Psi_{x}(\mathcal{M}_{|U}e^{-\frac{1}{x^q}})_{(0,0)}=0$ and $\Psi_{p_1}(\mathcal{M}e^{p_2-\alpha\circ p_1})$ has support included in $(0,\infty)$.

\vspace{0.5cm}

Let us now focus on $\Psi_{p_1}(\mathcal{M}e^{p_2-\alpha\circ p_1})_{(0,\infty)}$. In a neighbourhood $U$ of $(0,\infty)$ in $D\times\mathbb{P}^1$, there exist coordinates $(x,y)$ such that $p_1(x,y)=x$, $p_2(x,y)=\sfrac{1}{y}$ and $p_2-\alpha\circ p_1(x,y)=\sfrac{1}{y}-\alpha(x)$. To simplify notation, we will set $g=\sfrac{1}{y}-\alpha(x)$. Then we have:
$$
\Psi_{p_1}(\mathcal{M}e^{p_2-\alpha\circ p_1})_{(0,\infty)}=\Psi_{x}(\mathcal{M}_{|U}[\frac{1}{xy}]e^{g})_{(0,0)}.$$

We have to compute the Euler characteristic and the zeta function of the monodromy of $\DR~\Psi_{x}(\mathcal{M}_{|U}[\frac{1}{xy}]e^g)_{(0,0)}$.
We can not use directly the local computations because $g$ is not equivalent to a function of the type $\sfrac{1}{x^ky^l}$ at $(0,0)$. But we will be led to this situation with the help of a resolution. The choice of a good resolution is relevant for the end of the proof.

\subsubsection{Choice of a resolution}

\begin{lemma}\label{resolution}
There exists a resolution $\pi=(\pi_1,\pi_2):\mathbb{X}\to U$ with exceptional locus $E$ (actually it is a finite composition of the blow-up of the point $(0,0)\in U$ and some blow-up of points in its exceptional locus) which satisfies the following conditions:
\begin{enumerate}
\item Let $\widetilde{E}=\pi^{-1}(\{xy=0\})$. There exists a unique irreducible component $E_d$ of $E$ which intersects $\overline{\widetilde{E}\setminus E_d}$ in a unique point $P$. Moreover, in a neighbourhood of any point $Q\in E$, we can choose local coordinates $(u,v)$ on $\mathbb{X}$ such that:
\begin{enumerate}
\item If $Q\notin E_d$, 
\begin{itemize}
\item $\pi_1(u,v)=u^m$ and $g\circ\pi(u,v)\sim\sfrac{1}{u^k}$, with $m,k\geq 1$,

or 
\item $\pi_1(u,v)=u^mv^n$ and $g\circ\pi(u,v)\sim\sfrac{1}{u^kv^l}$, with $k,l,m\geq 1$, $n\geq 0$.
\end{itemize}
\item If $Q=P$,
\begin{itemize}
\item $u=0$ (resp. $v=0$) is an equation of $E_d$ (resp. of the other component of $E$),
\item $\pi_1(u,v)=uv^n$ and $g\circ\pi(u,v)\sim\sfrac{1}{v}$, with $n\geq 1$.
\end{itemize}
\item If $Q\in E_d\setminus\{P\}$, 
\begin{itemize}
\item $u=0$ is an equation of $E_d$,
\item $\pi_1(u,v)=u$ and $g\circ\pi(u,v)=a+v$, with $a\in\mathbb{C}$.
\end{itemize}
\end{enumerate}
\item The singular support of $\pi^{\ast}(\mathcal{M}_{|U}[\frac{1}{xy}])=\pi^{\ast}(\mathcal{M}_{|U})[\ast\widetilde{E}]$ has
at most normal crossings at any point of
$\overline{E\setminus E_d}$. 
\end{enumerate}
\end{lemma}

\newpage

\begin{figure}[h]
\includegraphics[height=3cm]{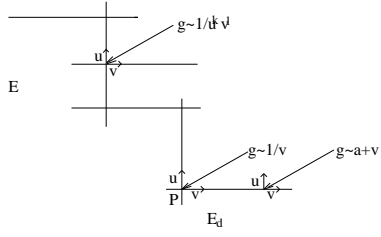}
\caption{Exceptional locus $E$}
\end{figure}

\begin{proof}
By reducing to the same denominator, we have $g=\sfrac{(x^q-y\beta)}{x^qy}$, with $\beta\in\mathbb{C}[x]$ and $\beta(0)\neq 0$.

Let us begin by proving the first point. We will perform $2q$ blow-up of points. The resolution tree we will obtain has just one branch and the component $E_d$ will be the exceptional locus of the last blow-up.

The first $q$ blow-up allow us to decrease the degree on $x$ in the numerator. We obtain a resolution $\pi=(\pi_1,\pi_2)$ composition of $q$ blow-up 
$$
\begin{array}{ccccccc}
X^{(q)}&\overset{\pi^{(q)}}{\to}&\ldots&\to&X^{(1)}&\overset{\pi^{(1)}}{\to}&U\\
\cup&&&&\cup&&\\
E_q&&&&E_1&&
\end{array},$$
with exceptional locus $E=E_1\cup\ldots\cup E_q$, such that $g\circ\pi$ can be written locally as in the case (a) except at a point $I\in E_q\setminus E_{q-1}$. In a neighbourhood of this point, we can choose local coordinates $(u,v)$ such that $\pi_1(u,v)=u$ and $g\circ\pi(u,v)=\sfrac{(u\gamma(u)-v\beta(u))}{u^qh(v)}$, with $h(v)=v+\sfrac{1}{\beta(0)}$ and $\gamma(u)=\sfrac{(\beta(0)-\beta(u))}{\beta(0)u}$. 

Now the $q-1$ next blow-up allow us to decrease the degree in $u$ in the denominator and the degree in $u$ of $\gamma(u)$. We obtain a new resolution $\pi=(\pi_1,\pi_2)$, composition of $2q-1$ blow-up, such that $g\circ\pi$ can be written locally as in the case (a) except at a point $J\in E_{2q-1}\setminus E_{2q-2}$. In some suitable local coordinates $(u,v)$ in a neighbourhood of $J$,  $\pi_1(u,v)=u$ and $g\circ\pi(u,v)=\sfrac{(cu-v\beta(u))}{uh(u,v)}$, with $c$ constant and $h(u,v)$ invertible.

At last, we blow-up $J$. The exceptional locus of this last blow up is $E_d$ and we can verify the local behaviour of $g$ on $E_d$. The details are left to the reader.

To prove the second point, it is sufficient to perform additional blow-up of points in $\overline{E\setminus E_d}$.
\end{proof}

We denote by $\widetilde{S_{\ell}}$ the strict transform of $S_{\ell}$.
\begin{corollary}\label{stricttransform}\mbox{}\par
$(1)$ Let $\ell\in\Lambda$. $\widetilde{S_{\ell}}\cap E_d\neq\emptyset$ if
and only if $\ell\in\Lambda_{\alpha}$.

$(2)$ The assumption $(\ast)$ is equivalent to:
$\forall \ell,\ell^{'}\in\Lambda_{\alpha},~\ell\neq\ell^{'}\Longrightarrow\widetilde{S_{\ell}}\cap E_d\neq\widetilde{S_{\ell^{'}}}\cap E_d$.
\end{corollary}

\begin{proof}
$(1)$ According to the choice of the resolution,
$\widetilde{S_{\ell}}\cap E_d\neq\emptyset$ if and only if
$\widetilde{S_{\ell}}\cap (E_d\setminus\{P\})\neq\emptyset$. Let $Q_{\ell}\in\widetilde{S_{\ell}}\cap (E_d\setminus\{P\})$ and  $(u,v)$ be some local coordinates in the
neighbourhood of $Q_{\ell}$ such that $f\circ\pi(u,v)=a+v$ and
$\pi_1(u,v)=u$, with $a\in\mathbb{C}$.

As the
intersection multiplicity of $S_{\ell}$ with $\{0\}\times\mathbb{P}^1$ is $1$ and
as $\pi_1(u,v)=u$, the intersection multiplicity of
$\widetilde{S_{\ell}}$ with $E_d$ is also $1$. Then:
\begin{displaymath}
\begin{split}
&\widetilde{S_{\ell}}\cap E_d\neq\emptyset,\\
&\Longleftrightarrow\exists\delta(s)\in\mathbb{C}\{s\}\text{ such that }a+v=\delta(u) \text{ is an equation of }\widetilde{S_{\ell}},\\
&\Longleftrightarrow\exists\delta(s)\in\mathbb{C}\{s\}\text{ such that }f\circ\pi(u,v)=\delta\circ\pi_1(u) \text{ is an equation of }\widetilde{S_{\ell}},\\
&\Longleftrightarrow\exists\delta(s)\in\mathbb{C}\{s\}\text{ such that }f(x,y)=\delta(x) \text{ is an equation of }S_{\ell},\\
&\Longleftrightarrow\exists\delta(s)\in\mathbb{C}\{s\}\text{ such that }\sfrac{1}{y}=\alpha(x)+\delta(x) \text{ is an equation of }S_{\ell},\\
&\Longleftrightarrow \ell\in\Lambda_{\alpha}.
\end{split}
\end{displaymath}

$(2)$ Here, we want to understand the behaviour of $S_{\ell}$ after the resolution $\pi$. We will keep the same notations as in the proof of Lemma \ref{resolution}.

Let $\ell\in\Lambda_{\alpha}$. Then $S_{\ell}$ has equation $h(x,y)=(\beta(x)+x^q\delta_{\ell}(x))y-x^q=0$. 

After the first $q$ blow-up, in the neighbourhood of the point $I\in E_q\setminus E_{q-1}$ with local coordinates $(u,v)$, 
$$h\circ\pi(u,v)=u^q\big((\beta(u)+u^q\delta_{\ell}(u))v+u(-\gamma(u)+u^{q-1}\sfrac{\delta_{\ell}(u)}{\beta(0)})\big).$$ 
After the next $(q-1)$ blow-up, in the neighbourhood of the point $J$, we can choose local coordinates $(u,v)$ such that 
$$h\circ\pi(u,v)=u^{2q-1}\big((\beta(u)+u^q\delta_{\ell}(u))v+u(-c+\sfrac{\delta_{\ell}(u)}{\beta(0)}+u\widetilde{\delta}(u))\big),$$
where $\widetilde{\delta}(u)\in\mathbb{C}\{u\}$ and $c$ is a constant which depends only on $\alpha$.

At last, we blow-up $J$ and obtain that in the chart $u=s$, $v=st$, 
$$h\circ\pi(s,t)=s^{2q}\big((\beta(s)+s^q\delta_{\ell}(s))t-c+\sfrac{\delta_{\ell}(s)}{\beta(0)}+s\widetilde{\delta}(s)\big).$$
Then, we see that the strict transform $\widetilde{S_{\ell}}$ of $S_{\ell}$ interesct $E_d$ at a point which depends only on $\alpha$ and $\delta_{\ell}(0)$. The assertion $(2)$ follows.
\end{proof}

\subsubsection{Use of the local computations}

There exists an isomorphism: 
$$\begin{array}{ll}
\DR~\Psi_{x}(\mathcal{M}_{|U}[\frac{1}{xy}]e^g)_{(0,0)}
&=R\Gamma(E,\DR~\Psi_{\pi_1}(\pi^{\ast}(\mathcal{M}_{|U}[\frac{1}{xy}])e^{g\circ\pi})),\\
&\text{(Same proof as for the isomorphism $(\ast\ast)$ p. \pageref{9}),}\\
&=R\Gamma(E,\DR~\Psi_{\pi_1}(\pi^{\ast}(\mathcal{M}_{|U})[\ast\widetilde{E}]e^{g\circ\pi})).
\end{array}$$

Let us first prove that the complex $\DR~\Psi_{\pi_1}(\pi^{\ast}(\mathcal{M}_{|U})[\ast\widetilde{E}]e^{g\circ\pi})_{|E}$ has support included in $E_d$. We have two cases to consider:

\begin{enumerate}
\item In a neighbourhood $V\subset\mathbb{X}$ of a smooth point $Q$ in $E\setminus E_d$, according to Lemma \ref{resolution}, there exist coordinates $(u,v)$ on $V$ such that $\pi_1(u,v)=u^m$, $g\circ\pi(u,v)\sim\sfrac{1}{u^k}$ and the singular support of $\pi^{\ast}(\mathcal{M}_{|U})[\ast\widetilde{E}]$ is included in $uv=0$. Then, according to Lemma \ref{B} 1, $$\Psi_{\pi_1}(\pi^{\ast}(\mathcal{M}_{|U})[\ast\widetilde{E}]e^{g\circ\pi})_Q=\Psi_{u^m}(\pi^{\ast}(\mathcal{M}_{|U})[\ast\widetilde{E}]e^{\frac{1}{u^k}})_{(0,0)}=0.$$

\item In a neighbourhood $V\subset\mathbb{X}$ of a normal crossing point $Q$ in $E\setminus E_d$, according to Lemma \ref{resolution}, there exist coordinates $(u,v)$ on $V$ such that $\pi_1(u,v)=u^mv^n$, $g\circ\pi(u,v)\sim\sfrac{1}{u^kv^l}$ and the singular support of $\pi^{\ast}(\mathcal{M}_{|U})[\ast\widetilde{E}]$ is included in $uv=0$. Then, according to Lemma \ref{B} 2, $$\Psi_{\pi_1}(\pi^{\ast}(\mathcal{M}_{|U})[\ast\widetilde{E}]e^{g\circ\pi})_Q=\Psi_{u^mv^n}(\pi^{\ast}(\mathcal{M}_{|U})[\ast\widetilde{E}]e^{\frac{1}{u^kv^l}})_{(0,0)}=0.$$
\end{enumerate}

Then $\DR~\Psi_{x}(\mathcal{M}_{|U}[\frac{1}{xy}]e^g)_{(0,0)}=R\Gamma(E_d,\DR~\Psi_{\pi_1}(\pi^{\ast}(\mathcal{M}_{|U})[\ast\widetilde{E}]e^{g\circ\pi}))$.

\vspace{0.5cm}

Now we want to examine the complex $\DR~\Psi_{\pi_1}(\pi^{\ast}(\mathcal{M}_{|U})[\ast\widetilde{E}]e^{g\circ\pi})$ locally on $E_d$. Let $T$ be a tubular neighbourhood of $E_d$. According to Corollary \ref{stricttransform} (1), the singular support $\widetilde{S_T}$ of $(\pi^{\ast}(\mathcal{M}_{|U})[\ast\widetilde{E}])_{|T}$ is $\big(E\cup(\cup_{\ell\in\Lambda_{\alpha}}\widetilde{S_\ell})\big)\cap T$.

\vspace{0.5cm}

\noindent
{\bf Computation of the Euler characteristic:} We denote by $Q_1,\ldots,Q_k$ the intersection point of $\cup_{\ell\in\Lambda_{\alpha}}\widetilde{S_\ell}$ with $E_d$.

Let $r$ be the rank of $\mathcal{M}$. As $\pi$ is an isomorphism out of $E$, the rank of $\pi^{\ast}(\mathcal{M}_{|U})$ is $r$ and the multiplicity of the conormal space $T_{\widetilde{S_\ell}}^{\ast}\mathbb{X}$ in the characteristic cycle of $\pi^{\ast}(\mathcal{M}_{|U})$ is $m_\ell$.

\begin{figure}[h]
\includegraphics[height=2.5cm]{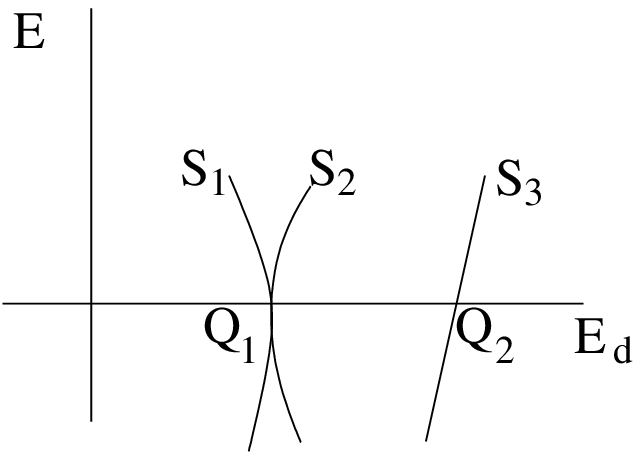}
\end{figure}

\begin{enumerate}
\item In a neighbourhood $V\subset\mathbb{X}$ of the point $P$, according to Lemma \ref{resolution}, there exist coordinates $(u,v)$ on $V$ such that $\pi_1(u,v)=uv^n$, $g\circ\pi(u,v)\sim\sfrac{1}{v}$ and the singular support of $\pi^{\ast}(\mathcal{M}_{|U})[\ast\widetilde{E}]$ is included in $uv=0$. Then,
$$\Psi_{\pi_1}(\pi^{\ast}(\mathcal{M}_{|U})[\ast\widetilde{E}]e^{g\circ\pi})_P=\Psi_{uv^n}(\pi^{\ast}(\mathcal{M}_{|U})[\frac{1}{uv}]e^{\frac{1}{v}})_{(0,0)}.$$ 

Applying Lemma \ref{B} 3., we conclude that the Euler characteristic of the complex $\DR~\Psi_{\pi_1}(\pi^{\ast}(\mathcal{M}_{|U})[\ast\widetilde{E}]e^{g\circ\pi})_P$ is equal to $-r$.

\item In a neighbourhood $V\subset\mathbb{X}$ of a point $Q_i$, $g\circ\pi$ is holomorphic and $\pi^{-1}_1(0)$ is smooth. Then we have
$\DR~\Psi_{\pi_1}(\pi^{\ast}(\mathcal{M}_{|U})[\ast\widetilde{E}]e^{g\circ\pi})_{Q_i}=\DR~\Psi_u(\pi^{\ast}(\mathcal{M}_{|U})[\frac{1}{u}])_{(0,0)}$
and, according to Lemma \ref{A}, the Euler characteristic of the complex $\DR~\Psi_{\pi_1}(\pi^{\ast}(\mathcal{M}_{|U})[\ast\widetilde{E}]e^{g\circ\pi})_{Q_i}$ is $r-\sum_{\widetilde{S_{\ell}}\cap E_d=\{Q_i\}}m_{\ell}$.

\item In a neighbourhood of a point $Q\in E_d\setminus\{P,Q_1,\ldots,Q_k\}$, as in the previous point, we prove that the Euler characteristic of the complex $\DR~\Psi_{\pi_1}(\pi^{\ast}(\mathcal{M}_{|U})[\ast\widetilde{E}]e^{g\circ\pi})_Q$ is $r$ (because $\{\ell\in\Lambda~|~\widetilde{S_\ell}\cap E_d=\{Q\}\}=\emptyset$).
\end{enumerate}

Using the Mayer-Vietoris Theorem, we deduce that the Euler characteristic of the complex $R\Gamma(E_d,\DR~\Psi_{\pi_1}(\pi^{\ast}(\mathcal{M}_{|U})[\ast\widetilde{E}]e^{g\circ\pi}))$ is equal to

\begin{displaymath}
\begin{split}
-r+\sum_{i=1}^k\bigg(r-\sum_{\widetilde{S_{\ell}}\cap E_d=\{Q_i\}}m_{\ell}\bigg)+\chi(E_d\setminus\{P,Q_1,\ldots,Q_k\}).r&=
-\sum_{\widetilde{S_{\ell}}\cap E_d\neq\emptyset}m_{\ell}\\
&=-\sum_{\ell\in\Lambda_{\alpha}}m_{\ell}.
\end{split}
\end{displaymath}

{\bf Computation of the zeta function of the monodromy:} Assume that $(\ast)$ is fulfilled. The fact that for any $\ell\in\Lambda$, $p_\ell=1$ and the assumption $(\ast)$ imply that $\widetilde{S_T}$ is a divisor with normal crossing (cf. Corollary \ref{stricttransform} (2)).

We denote by $T^{\ast}$ the complement of $\widetilde{S_T}$ in $T$ and by $Q_{\ell}$ the intersection point of one $\widetilde{S_\ell}$ with $E_d$, $\ell\in\Lambda_{\alpha}$.

\begin{figure}[h]
\includegraphics[height=2.5cm]{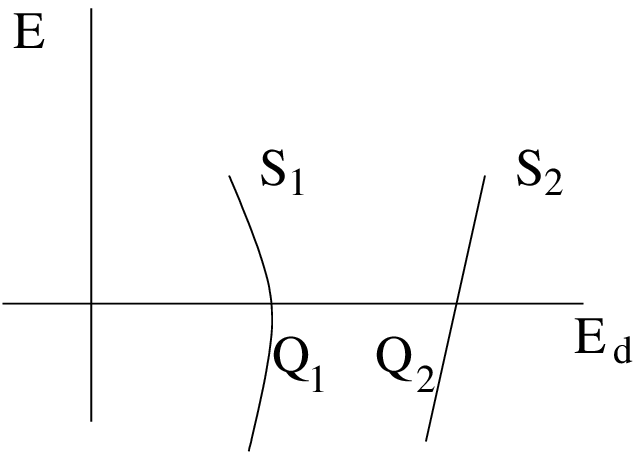}
\end{figure}

\begin{itemize}
\item[$\bullet$] $(\DR~\pi^{\ast}(\mathcal{M}_{|U})[\ast\widetilde{E}])_{|T^{\ast}}$ is a local system. We can define the monodromy endomorphism around $E_d$, $T_P:(\DR~\pi^{\ast}(\mathcal{M}_{|U})[\ast\widetilde{E}])_{P}\to (\DR~\pi^{\ast}(\mathcal{M}_{|U})[\ast\widetilde{E}])_{P}$, for any $P\in T^{\ast}$. As $T^{\ast}$ is connected and $\widetilde{S_T}$ has normal crossings, the characteristic polynomial of $T_P$ does not depend on $P\in T^{\ast}$. We denote it by $\zeta_r$. 
\item[$\bullet$] To any $\ell\in\Lambda_{\alpha}$, we associate the characteristic polynomial $\zeta_\ell$ of the monodromy of the local system $\phi_{\widetilde{h_\ell}}(\DR~\pi^{\ast}(\mathcal{M}_{|U}))_{|\widetilde{S_\ell}\setminus\{(0,0)\}}$ around $(0,0)$, where $\widetilde{h_\ell}$ is an equation of $\widetilde{S_\ell}$. If $h_\ell$ is an equation of $S_\ell$, as $\pi$ is an isomorphism out of $E$, it coincides with the characteristic polynomial of the monodromy of $\phi_{h_\ell}(\DR~\mathcal{M})_{|S_\ell\setminus\{(0,0)\}}$ around $(0,0)$.
\end{itemize}

Using the local behaviour of $\pi_1$ and $g\circ\pi$ on $E_d$ given by Lemma \ref{resolution}, we have:

\begin{enumerate}
\item In a neighbourhood of the point $P$, according to Lemma \ref{B} 3., the zeta function of the monodromy of $\DR~\Psi_{\pi_1}(\pi^{\ast}(\mathcal{M}_{|U})[\ast\widetilde{E}]e^{g\circ\pi})_P$ is equal to $\zeta_r^{-1}$.
\item In a neighbourhood of a point $Q_\ell$, $\ell\in\Lambda_{\alpha}$, according to Lemma \ref{C}, the zeta function of the monodromy of $\DR~\Psi_{\pi_1}(\pi^{\ast}(\mathcal{M}_{|U})[\ast\widetilde{E}]e^{g\circ\pi})_{Q_\ell}$ is equal to $\zeta_r\zeta_\ell^{-1}$.
\item In a neighbourhood of a point $Q\in E_d\setminus\{P,Q_\ell,\ell\in\Lambda_{\alpha}\}$, the zeta function of the monodromy of $\DR~\Psi_{\pi_1}(\pi^{\ast}(\mathcal{M}_{|U})[\ast\widetilde{E}]e^{g\circ\pi})_Q$ is equal to $\zeta_r$.
\end{enumerate}

Then, using the Mayer-Vietoris Theorem, we deduce that the zeta function of the monodromy of $R\Gamma(E_d,\DR~\Psi_{\pi_1}(\pi^{\ast}(\mathcal{M}_{|U})[\ast\widetilde{E}]e^{g\circ\pi}))$ is equal to
$$\zeta_r^{-1}.\prod_{\ell\in\Lambda_{\alpha}}(\zeta_r\zeta_\ell^{-1}).\zeta_r^{\chi(E_d\setminus\{P,Q_\ell,\ell\in\Lambda_{\alpha}\})}=\prod_{\ell\in\Lambda_{\alpha}}\zeta_\ell^{-1}.$$

\section{Realization of  $\mathbb{C}[[t]]\langle\partial_t\rangle$-modules}\label{realization}
Let $\widehat{\mathcal{N}}$ be a formal holonomic $\mathbb{C}[[t]]\langle\partial_t\rangle$-module.

\begin{theorem}\label{realisation}
There exist a ramification $\rho:\tau\to t=\tau^p$ and a regular holonomic $\mathcal{D}_{D\times\mathbb{P}^1}$-module $\mathcal{M}$ such that $\rho^{\ast}(\widehat{\mathcal{N}})$ is isomorphic to the formalization of $\rho^{\ast}(\mathcal{H}^0p_{1+}(\mathcal{M}e^{p_2}))_0$.
\end{theorem}

\begin{proof}
We recall that we identify $\mathcal{D}_{D,0}$ with $\mathbb{C}\{t\}\langle\partial_t\rangle$, by choosing a coordinate $t$ on $D$. We begin by stating two lemmas.
 
\begin{lemma}\label{BB} 
\begin{enumerate}
\item $p_{1+}(\mathcal{O}_{D\times\mathbb{P}^1}e^{p_2})=\mathcal{O}_{D}$.
\item $p_{1+}(\mathcal{O}_{D\times\mathbb{P}^1}[t^{-1}]e^{p_2})=\mathcal{O}_{D}[t^{-1}]$.
\item Let $(p,\alpha)\in\mathbb{N}^{\ast}\times(\tau^{-1}.\mathbb{C}[\tau^{-1}])$. Let $\rho:D^{'}\to D$, defined by $\rho(\tau)=t=\tau^p$. We denote by $Z$ the curve in $D\times\mathbb{P}^1$ parametrized by $\tau\to(\tau^p,\alpha)$ and by $\widetilde{Z}$ the union of $Z$, $\{0\}\times\mathbb{P}^1$ and $D\times\{\infty\}$.
$$\rho^{\ast}p_{1+}(\mathcal{O}_{D\times\mathbb{P}^1}[\ast \widetilde{Z}]e^{p_2})=\underset{\xi^p=1}{\oplus}\mathcal{O}_{D^{'}}[\tau^{-1}]e^{\alpha(\xi\tau)}\oplus \mathcal{O}_{D^{'}}[\tau^{-1}].$$ 
\end{enumerate}
\end{lemma}

\begin{proof}
We just prove the third point. The first two points can be proved similary.

Let $p_1^{'}:D^{'}\times\mathbb{P}^1\to D^{'}$ and $p_2^{'}:D^{'}\times\mathbb{P}^1\to\mathbb{P}^1$ be the canonical projections and $\rho^{'}=(\rho,id):D^{'}\times\mathbb{P}^1\to D\times\mathbb{P}^1$.

First note that $\rho^{\ast}p_{1+}(\mathcal{O}_{D\times\mathbb{P}^1}[\ast \widetilde{Z}]e^{p_2})=p^{'}_{1+}(\mathcal{O}_{D^{'}\times\mathbb{P}^1}[\ast \rho^{'-1}(\widetilde{Z})]e^{p_2^{'}})$ (cf. Lemma \ref{comm}).

Denote by $\mathcal{K}^{\bullet}:=~^pDR_{D^{'}\times\mathbb{P}^1/D^{'}}(\mathcal{O}_{D^{'}\times\mathbb{P}^1}[\ast \rho^{'-1}(\widetilde{Z})]e^{p_2^{'}})$ the relative de Rham complex of $p_1$ concentrated in negative degrees. It consists in the complex $$\mathcal{O}_{D^{'}\times\mathbb{P}^1}[\ast \rho^{'-1}(\widetilde{Z})]\overset{d}{\to}\mathcal{O}_{D^{'}\times\mathbb{P}^1}[\ast \rho^{'-1}(\widetilde{Z})],$$
with $d(h)=(\partial_{\mathbb{P}^1}(h)+h\partial_{\mathbb{P}^1}(p_2^{'}))$.

By definition, $p_{1+}^{'}(\mathcal{O}_{D^{'}\times\mathbb{P}^1}[\ast \rho^{'-1}(\widetilde{Z})]e^{p_2^{'}})=\mathbb{R}p_{1\ast}^{'}(\mathcal{K}^{\bullet})$. We consider the spectral sequence $\mathcal{H}^{i-k}(R^kp_{1\ast}^{'}(\mathcal{K}^{\bullet}),d)\Longrightarrow \mathbb{R}^ip_{1\ast}^{'}(\mathcal{K}^{\bullet})$, where $(R^kp_{1\ast}^{'}(\mathcal{K}^{\bullet}),d)$ denotes the complex $R^kp_{1\ast}^{'}(\mathcal{O}_{D^{'}\times\mathbb{P}^1}[\ast \rho^{'-1}(\widetilde{Z})])\overset{d}{\to}R^kp_{1\ast}^{'}(\mathcal{O}_{D^{'}\times\mathbb{P}^1}[\ast \rho^{'-1}(\widetilde{Z})])$.

We claim that $R^kp_{1\ast}^{'}(\mathcal{O}_{D^{'}\times\mathbb{P}^1}[\ast \rho^{'-1}(\widetilde{Z})])=0$, for $k\neq 0$. Indeed, let $j:(D^{'}\times\mathbb{P}^1)\setminus\rho^{'-1}(\widetilde{Z})\hookrightarrow D^{'}\times\mathbb{P}^1$ and $P\in D^{'}$:
$$\begin{array}{ll}
R^kp_{1\ast}^{'}&(\mathcal{O}_{D^{'}\times\mathbb{P}^1}[\ast \rho^{'-1}(\widetilde{Z})])_P=\\
&=R^kp_{1\ast}^{'}Rj_{\ast}j^{-1}(\mathcal{O}_{D^{'}\times\mathbb{P}^1})_P,\\
&=R^k(p_1^{'}\circ j)_{\ast}j^{-1}(\mathcal{O}_{D^{'}\times\mathbb{P}^1})_P,\\
&=\varinjlim_{U\ni P open}R\Gamma(j^{-1}(U\times\mathbb{P}^1),j^{-1}(\mathcal{O}_{D^{'}\times\mathbb{P}^1})).\\
\end{array}$$

As $j^{-1}(U\times\mathbb{P}^1)=(U\times\mathbb{C})\setminus\rho^{'-1}(Z\cup\{0\}\times\mathbb{C})$ is Stein, we conclude that  $R^kp_{1\ast}^{'}(\mathcal{O}_{D^{'}\times\mathbb{P}^1}[\ast \rho^{'-1}(\widetilde{Z})])=0$, for $k\neq 0$.

Now $(R^0p_{1\ast}^{'}(\mathcal{K}^{\bullet}),d)$ is the complex:
$$\mathcal{O}_{D^{'}}[y]\bigg[\frac{1}{\tau\underset{\xi^p=1}{\prod}(y-\alpha(\xi\tau))}\bigg]\overset{d}{\to} \mathcal{O}_{D^{'}}[y]\bigg[\frac{1}{\tau\underset{\xi^p=1}{\prod}(y-\alpha(\xi\tau))}\bigg],$$
with $d(h)=(\partial_y(h)+h)$.

\begin{itemize}
\item Computation of $\mathcal{H}^{-1}(R^0p_{1\ast}^{'}(\mathcal{K}^{\bullet}),d)$:\\
Let $h\in\mathcal{O}_{D^{'}}[y]\bigg[\frac{1}{\tau\underset{\xi^p=1}{\prod}(y-\alpha(\xi\tau))}\bigg]$ such that $d(h)=0$. Thus $h=Ce^{-y}$, with $C\in\mathcal{O}_{D^{'}}[\frac{1}{\tau}]$ and then $h=0$. Thus, $\mathcal{H}^{-1}(R^0p_{1\ast}^{'}(\mathcal{K}^{\bullet}),d)=0$.
\item Computation of $\mathcal{H}:=\mathcal{H}^{0}(R^0p_{1\ast}^{'}(\mathcal{K}^{\bullet}),d)$:\\
We note first that $\mathcal{H}$ is generated over $\mathcal{O}_{D^{'}}[\frac{1}{\tau}]$ by \\$\{y^m,\frac{y^m}{(y-\alpha(\xi\tau))^n};~m\geq 0,~n\geq 1 \text{ and }\xi^p=1\}$.

Moreover, if $m,n>0$, we have:
$\frac{y^m}{(y-\alpha(\xi\tau))^n}=\frac{y^{m-1}}{(y-\alpha(\xi\tau))^{n-1}}+\frac{y^{m-1}\alpha(\xi\tau)}{(y-\alpha(\xi\tau))^n}$ and $y^m\equiv -my^{m-1}$.

Then $\mathcal{H}$ is generated over $\mathcal{O}_{D^{'}}[\frac{1}{\tau}]$ by $\{1,\frac{1}{(y-\alpha(\xi\tau))^n};~n\geq 1 \text{ and }\xi^p=1\}$. But we have $\frac{1}{(y-\alpha(\xi\tau))^n}\equiv\frac{n}{(y-\alpha(\xi\tau))^{n+1}}$. Then, $\mathcal{H}$ is generated over $\mathcal{O}_{D^{'}}[\frac{1}{\tau}]$ by $\{1,\frac{1}{(y-\alpha(\xi\tau))};~\xi^p=1\}$.

Then we have  $\mathcal{H}=\underset{\xi^p=1}{\oplus}\mathcal{O}_{D^{'}}[\tau^{-1}].\frac{1}{(y-\alpha(\xi\tau))}\oplus \mathcal{O}_{D^{'}}[\tau^{-1}]$. 

As $
\partial_{\tau}(\frac{1}{y-\alpha(\xi\tau)})=\partial_{\tau}(\alpha(\xi\tau)).\frac{1}{(y-\alpha(\xi\tau))^2}\equiv\partial_{\tau}(\alpha(\xi\tau)).\frac{1}{y-\alpha(\xi\tau)}$, we have an isomorphism of $\mathcal{D}_{D^{'}}$-modules
$$\mathcal{O}_{D^{'}}[\tau^{-1}].\frac{1}{(y-\alpha(\xi\tau))}=\mathcal{O}_{D^{'}}[\tau^{-1}]e^{\alpha(\xi\tau)}.$$

Lemma \ref{BB} follows.
\end{itemize}

\end{proof}

\begin{lemma}\label{AA}
Let $\mathcal{N}$ be a holonomic $\mathcal{D}_{D}$-module and $\mathcal{M}$ be a holonomic $\mathcal{D}_{D\times\mathbb{P}^1}$-module. Suppose that $\mathcal{M}$ is $\mathcal{O}_{D\times\mathbb{P}^1}$-flat. We have an isomorphism: $$p_{1+}(\mathcal{M}\otimes_{\mathcal{O}_{D\times\mathbb{P}^1}}p_1^{\ast}(\mathcal{N}))=p_{1+}(\mathcal{M})\otimes_{\mathcal{O}_{D}}^{\mathbb{L}}\mathcal{N}.$$
\end{lemma}

\begin{proof}
First, we remark that $\mathcal{M}\otimes_{\mathcal{O}_{D\times\mathbb{P}^1}}p_1^{\ast}\mathcal{N}=\mathcal{M}\otimes_{p_1^{-1}\mathcal{O}_D}p_1^{-1}\mathcal{N}$, with the structure of $\mathcal{D}_{D\times\mathbb{P}^1}$-module given by $\xi(m\otimes n)=\xi(m)\otimes n+\xi(p_1)m\otimes \partial_t(n)$, where $\xi$ is a local section of the sheaf of vector fields on $D\times\mathbb{P}^1$.

Denote by $~^p\DR_{D\times\mathbb{P}^1/D}$ (resp. $\Omega^{\bullet+1}_{D\times\mathbb{P}^1/D}$) the relative de Rham functor of $p_1$ (resp. the complex of relative forms of $p_1$). The differential $d$ of the complex $~^p\DR_{D\times\mathbb{P}^1/D}(\mathcal{M}\otimes_{p_1^{-1}\mathcal{O}_D}p_1^{-1}\mathcal{N})=\Omega^{\bullet+1}_{D\times\mathbb{P}^1/D}\otimes_{\mathcal{O}_{D\times\mathbb{P}^1}}(\mathcal{M}\otimes_{p_1^{-1}\mathcal{O}_D}p_1^{-1}\mathcal{N})$ is defined by:
$$\begin{array}{lll}
d(w\otimes(m\otimes n))&=&dw\otimes(m\otimes n)+dp_2\wedge w\otimes\partial_{\mathbb{P}^1}(m\otimes n),\\
&=&dw\otimes(m\otimes n)+dp_2\wedge w\otimes(\partial_{\mathbb{P}^1}(m)\otimes n),
\end{array}$$
and the action of $p_1^{-1}\mathcal{D}_D$ is given by:
$$\partial_t(w\otimes(m\otimes n))=w\otimes(\partial_tm\otimes n)+w\otimes(m\otimes\partial_tn).$$

Then, we have an isomorphism of $p_1^{-1}\mathcal{D}_D$-modules:
$$
~^p\DR_{D\times\mathbb{P}^1/D}(\mathcal{M}\otimes_{p_1^{-1}\mathcal{O}_D}p_1^{-1}\mathcal{N})=~^p\DR_{D\times\mathbb{P}^1/D}(\mathcal{M})\otimes_{p_1^{-1}\mathcal{O}_D}p_1^{-1}\mathcal{N}.$$

We conclude with:
$$\begin{array}{ll}
p_{1+}(\mathcal{M}&\otimes^{\mathbb{L}}_{\mathcal{O}_{D\times\mathbb{P}^1}}p_1^{\ast}(\mathcal{N}))=\\
&=\mathbb{R}p_{1\ast}(~^p\DR_{D\times\mathbb{P}^{1}/D}(\mathcal{M}\otimes_{\mathcal{O}_{D\times\mathbb{P}^1}}p_1^{\ast}(\mathcal{N}))),\\
&=\mathbb{R}p_{1\ast}(~^p\DR_{D\times\mathbb{P}^{1}/D}(\mathcal{M})\otimes_{p_1^{-1}\mathcal{O}_D}p_1^{-1}(\mathcal{N})),\\
&=\mathbb{R}p_{1\ast}(~^p\DR_{D\times\mathbb{P}^{1}/D}(\mathcal{M})\otimes^{\mathbb{L}}_{p_1^{-1}\mathcal{O}_D}p_1^{-1}(\mathcal{N})),\\
&(~^p\DR_{D\times\mathbb{P}^{1}/D}(\mathcal{M})\text{ is }p_1^{-1}\mathcal{O}_D\text{-flat),}\\
&=\mathbb{R}p_{1\ast}(~^p\DR_{D\times\mathbb{P}^{1}/D}(\mathcal{M})\otimes^{\mathbb{L}}_{p_1^{-1}\mathcal{D}_D}(p_1^{-1}\mathcal{D}_D\otimes_{p_1^{-1}\mathcal{O}_D}p_1^{-1}(\mathcal{N}))),\\
&=\mathbb{R}p_{1\ast}(~^p\DR_{D\times\mathbb{P}^{1}/D}(\mathcal{M}))\otimes^{\mathbb{L}}_{\mathcal{D}_D}(\mathcal{D}_D\otimes_{\mathcal{O}_D}\mathcal{N}),\\
&\text{(Proposition 2.6.6 in \cite{KaSc}),}\\
&p_{1+}(\mathcal{M})\otimes^{\mathbb{L}}_{\mathcal{O}_D}\mathcal{N}.
\end{array}$$
\end{proof}

According to the structure of a formal $\mathbb{C}[[t]]\langle\partial_t\rangle$-module (cf. \cite{Ma}), we are led to prove the theorem in the cases of regular modules and purely irregular modules.
\begin{enumerate}
\item $\underline{\text{Regular case:}}$ Suppose $\widehat{\mathcal{N}}$ regular.
There exist a small disc $D$ in $\mathbb{C}$ centered at the origin and a regular holonomic $\mathcal{D}_D$-module $\mathcal{N}$, such that the formalization of $\mathcal{N}_0$ is $\widehat{\mathcal{N}}$ (cf. Theorem 5.3 p. 38 in \cite{Ma}).
Let $\mathcal{M}$ be the regular holonomic $\mathcal{D}_{D\times\mathbb{P}^1}$-module $p_1^{\ast}(\mathcal{N})$. We have:
$$\begin{array}{ll}
p_{1+}(\mathcal{M}e^{p_2})&=p_{1+}(\mathcal{O}_{D\times\mathbb{P}^1}e^{p_2}\otimes_{\mathcal{O}_{D\times\mathbb{P}^1}}p_1^{\ast}(\mathcal{N})),\\
&=p_{1+}(\mathcal{O}_{D\times\mathbb{P}^1}e^{p_2})\otimes_{\mathcal{O}_D}^{\mathbb{L}}\mathcal{N},\text{ (Lemma \ref{AA}),}\\
&=\mathcal{O}_D\otimes_{\mathcal{O}_D}^{\mathbb{L}}\mathcal{N},\text{ (Lemma \ref{BB}),}\\
&=\mathcal{N}.
\end{array}$$

\item $\underline{\text{Purely irregular case:}}$ Let $(p,\alpha)\in \mathbb{N}^{\ast}\times(\tau^{-1}.\mathbb{C}[\tau^{-1}])$ and $\widehat{R_{\alpha}}$ be a regular holonomic $\mathbb{C}[[\tau]]\langle\partial_{\tau}\rangle$-module, such that $\widehat{R_{\alpha}}[\tau^{-1}]=\widehat{R_{\alpha}}$. Let $\rho:\tau\to t=\tau^p$ and suppose $\rho^{\ast}(\widehat{\mathcal{N}})=\underset{\xi^p=1}{\oplus}\widehat{R_{\alpha}}e^{\alpha(\xi\tau)}$.

There exist a small disc $D^{'}$ in $\mathbb{C}$ centered at the origin and a regular holonomic $\mathcal{D}_{D^{'}}$-module $\widetilde{R_{\alpha}}$ such that the formalization of $(\widetilde{R_{\alpha}})_0$ is $\widehat{R_{\alpha}}$ (cf. Theorem 5.3 p. 38 in \cite{Ma}).But there exists a regular holonomic $\mathcal{D}_{D}$-module $R_{\alpha}$ such that $\widetilde{R_{\alpha}}=\rho^{\ast}(R_{\alpha})$, where $\rho:D^{'}\to D$, $\rho(\tau)=t=\tau^p$. Indeed, as $\rho^{\ast}$ is exact and $\widetilde{R_{\alpha}}$ is regular holonomic with $\widetilde{R_{\alpha}}[\tau^{-1}]=\widetilde{R_{\alpha}}$, we can reduce the proof to the case where 
$\widetilde{R_{\alpha}}=\sfrac{\mathcal{D}_{D^{'}}}{\mathcal{D}_{D^{'}}((\tau\partial_{\tau}-\alpha)^r)}$, with $r\in\mathbb{N}^{\ast}$ and $\alpha\in\mathbb{C}\setminus(-\mathbb{N})$ (cf. Corollary 19 in \cite{BrMai}). Then, we consider $R_{\alpha}=\sfrac{\mathcal{D}_D}{\mathcal{D}_D((t\partial_t-\frac{\alpha}{p})^r)}$.

Using notations of Lemma \ref{BB}, we denote by $\mathcal{M}$ the regular holonomic  $\mathcal{D}_{D\times\mathbb{P}^1}$-module $\sfrac{p_1^{\ast}(R_{\alpha})[\ast \widetilde{Z}]}{p_1^{\ast}(R_{\alpha})[t^{-1}]}$. 
According to Lemmas \ref{AA} and \ref{BB}, we have two isomorphisms:
$$\begin{array}{l}
\rho^{\ast}p_{1+}(p_1^{\ast}(R_{\alpha})[\ast \widetilde{Z}]e^{p_2})=\underset{\xi^p=1}{\oplus}\widetilde{R_{\alpha}}e^{\alpha(\xi\tau)}\oplus \widetilde{R_{\alpha}}\\
\rho^{\ast}p_{1+}(p_1^{\ast}(R_{\alpha})[t^{-1}]e^{p_2})=\widetilde{R_{\alpha}}.
\end{array}$$

Then the formalization of $\rho^{\ast}p_{1+}(\mathcal{M}e^{p_2})$ is $\underset{\xi^p=1}{\oplus}\widehat{R_{\alpha}}e^{\alpha(\xi\tau)}$.
\end{enumerate}

\end{proof}

\section{Direct image of holonomic $\mathcal{D}$-module of exponential type: the algebraic case}\label{algebraic}

In this section, we prove Theorem \ref{2}.
We consider $\mathcal{N}^k:=\mathcal{H}^kj_+f_+(\mathcal{M}e^g)$, the extension of $\mathcal{H}^kf_+(\mathcal{M}e^g)$ at infinity, where $j:\mathbb{A}^1\hookrightarrow\mathbb{A}^1\cup\{\infty\}\simeq\mathbb{P}^1$. Let $c\in\mathbb{P}^1$. The proof will be divided in two steps.
\begin{enumerate}
\item We consider the following diagrams
$$\xymatrix{\mathbb{C}^2\ar@{^(->}[r]^-i\ar[d]^-{\pi_1}&\mathbb{P}^1\times\mathbb{P}^1\ar[d]^-{p_1}&B_c\times\mathbb{P}^1\ar@{_(->}[l]_-{i_c}\ar[d]^{p_1}
&\mathbb{C}^2\ar@{^(->}[r]^-i\ar[d]^-{\pi_2}&\mathbb{P}^1\times\mathbb{P}^1\ar[d]^-{p_2}&B_c\times\mathbb{P}^1\ar@{_(->}[l]_-{i_c}\ar[d]^{p_2}\\
\mathbb{C}\ar@{^(->}[r]^-j&\mathbb{P}^1&B_c\ar@{_(->}[l]_-{j_c},&\mathbb{C}\ar@{^(->}[r]^-j&\mathbb{P}^1&\mathbb{P}^1\ar@{=}[l],}$$
where $\pi_1$, $\pi_2$ (resp. $p_1$, $p_2$) are the two canonical
projections. 

Then we have:
$$\begin{array}{lll}
j_c^{\ast}(j_+f_+(\mathcal{M}e^g))^{\an}&=&j_c^{\ast}(j_+\pi_{1+}(f,g)_+(\mathcal{M}e^g))^{\an},\\
&&\text{(Prop. 6.4 in \cite{Bo} and }f=\pi_1\circ(f,g)),\\
&=&j_c^{\ast}(p_{1+}i_+(f,g)_+(\mathcal{M}e^g))^{\an},\\
&&\text{(Prop. 6.4 in \cite{Bo} and }j\circ\pi_1=p_1\circ i),\\
&=&j_c^{\ast}p_{1+}(i_+(f,g)_+(\mathcal{M}e^g))^{\an},\\
&&\text{(Prop 8.2.2 p. 179 in \cite{Me})},\\
&=&p_{1+}i_c^{\ast}((i_+(f,g)_+(\mathcal{M}))e^{p_2})^{\an},\\
&&(i_c\text{ is an open inclusion}),\\
&=&p_{1+}(\mathcal{P}_ce^{p_2}).
\end{array}$$

As the formalization of the germ $\mathcal{N}_c^k$ and of its analytization are equal, $\mathcal{N}_c^k$ and $\mathcal{H}^kp_{1+}(\mathcal{P}_ce^{p_2})_c$ have the same formal irregular part.

\item Then we claim that $\mathcal{H}^kp_{1+}(\mathcal{P}_ce^{p_2})_c$ and $\mathcal{H}^0p_{1+}(\mathcal{H}^k\mathcal{P}_ce^{p_2})_c$ have the same formal irregular part. According to the formal decomposition Theorem (cf. Theorem 1.2 p. 43 and Theorem 2.3 p. 51 in \cite{Ma}), there exists a ramification $\rho:D^{'}_c\to D_c$, $\rho(\tau)=\tau^p=t$, such that the modules $\rho^{\ast}(\mathcal{H}^kp_{1+}(\mathcal{P}_ce^{p_2}))_c$ and $\rho^{\ast}(\mathcal{H}^0p_{1+}(\mathcal{H}^k\mathcal{P}_ce^{p_2}))_c$ are isomorphic to the finite direct sum of some formal modules of exponential type. Then, it is sufficient to prove that, given an element $\alpha\in\tau^{-1}.\mathbb{C}[\tau^{-1}]$, the regular part of  $(\rho^{\ast}(\mathcal{H}^kp_{1+}(\mathcal{P}_ce^{p_2}))_c)e^{-\alpha}$ and  $(\rho^{\ast}(\mathcal{H}^0p_{1+}(\mathcal{H}^k\mathcal{P}_ce^{p_2}))_c)e^{-\alpha}$ are isomorphic. Thus, we are led to prove that the specialization at $c$ of  $(\rho^{\ast}(\mathcal{H}^kp_{1+}(\mathcal{P}_ce^{p_2})))e^{-\alpha}$ and  $(\rho^{\ast}(\mathcal{H}^0p_{1+}(\mathcal{H}^k\mathcal{P}_ce^{p_2})))e^{-\alpha}$ are isomorphic (cf. Example 5.2.1 in \cite{LaMa}).

Let $\rho^{'}=(\rho,id):D^{'}_c\times\mathbb{P}^1\to  D_c\times\mathbb{P}^1$. We denote by $p_1^{'}:D_c\times\mathbb{P}^1\to D_c$ and $p_2^{'}:D_c\times\mathbb{P}^1\to\mathbb{P}^1$ the canonical projection. Let $\overline{p}_1^{'}:\{c\}\times\mathbb{P}^1\to\{c\}$. We have:
$$\begin{array}{ll}
\spe_c((\rho^{\ast}&(\mathcal{H}^kp_{1+}(\mathcal{P}_ce^{p_2})))e^{-\alpha})_c\\
&=\spe_c(\mathcal{H}^kp_{1+}^{'}(\rho^{'\ast}(\mathcal{P}_c)e^{p_2^{'}})e^{-\alpha})_c,\\
&\text{(Lemma \ref{comm} extended to complexes of $\mathcal{D}$-modules),}\\
&=\spe_c(\mathcal{H}^kp_{1+}^{'}(\rho^{'\ast}(\mathcal{P}_c)e^{p_2^{'}-\alpha\circ p_1^{'}}))_c,\\
&=\mathcal{H}^k\overline{p}_{1+}^{'}\spe_{\{c\}\times\mathbb{P}^1}(\rho^{'\ast}(\mathcal{P}_c)e^{p_2^{'}-\alpha\circ p_1^{'}})_c,\\
&\text{(Theorem 9.4.1 in \cite{LaMa}),}\\
&=R^k\Gamma(\{c\}\times\mathbb{P}^1,\DR~\spe_{\{c\}\times\mathbb{P}^1}(\rho^{'\ast}(\mathcal{P}_c)e^{p_2^{'}-\alpha\circ p_1^{'}})[+1]),\\
&\text{(p.5 in \cite{Ma3}).}
\end{array}$$

According to Proposition \ref{11} extended to complexes of $\mathcal{D}$-modules, the support of  $\spe_{\{c\}\times\mathbb{P}^1}(\rho^{'\ast}(\mathcal{P}_c)e^{p_2^{'}-\alpha\circ p_1^{'}})$ is included in $(c,\infty)$. Then,

$$\begin{array}{ll}
\spe_c&((\rho^{\ast}(\mathcal{H}^kp_{1+}(\mathcal{P}_ce^{p_2})))e^{-\alpha})_c=\\
&=(\mathcal{H}^k\DR~\spe_{\{c\}\times\mathbb{P}^1}(\rho^{'\ast}(\mathcal{P}_c)e^{p_2^{'}-\alpha\circ p_1^{'}})[+1])_{(c,\infty)},\\
&=(\DR~\mathcal{H}^k\spe_{\{c\}\times\mathbb{P}^1}(\rho^{'\ast}(\mathcal{P}_c)e^{p_2^{'}-\alpha\circ p_1^{'}})[+1])_{(c,\infty)},\\
&=R^0\Gamma(\{c\}\times\mathbb{P}^1,\DR~\mathcal{H}^k\spe_{\{c\}\times\mathbb{P}^1}(\rho^{'\ast}(\mathcal{P}_c)e^{p_2^{'}-\alpha\circ p_1^{'}})[+1]),\\
&=R^0\Gamma(\{c\}\times\mathbb{P}^1,\DR~\spe_{\{c\}\times\mathbb{P}^1}(\rho^{'\ast}(\mathcal{H}^k\mathcal{P}_c)e^{p_2^{'}-\alpha\circ p_1^{'}})[+1]),\\
&=\mathcal{H}^0\overline{p}_{1+}^{'}\spe_{\{c\}\times\mathbb{P}^1}(\rho^{'\ast}(\mathcal{H}^k\mathcal{P}_c)e^{p_2^{'}-\alpha\circ p_1^{'}})_c,\\
&\text{(p.5 in \cite{Ma3}),}\\
&=\spe_c(\mathcal{H}^0p_{1+}^{'}(\rho^{'\ast}(\mathcal{H}^k\mathcal{P}_c)e^{p_2^{'}-\alpha\circ p_1^{'}}))_c,\\
&\text{(Theorem 4.8.1 p. 226 in \cite{MeSa}),}\\
&=\spe_c(\mathcal{H}^0p_{1+}^{'}(\rho^{'\ast}(\mathcal{H}^k\mathcal{P}_c)e^{p_2^{'}})e^{-\alpha})_c,\\
&=\spe_c((\rho^{\ast}(\mathcal{H}^0p_{1+}(\mathcal{H}^k\mathcal{P}_ce^{p_2})))e^{-\alpha})_c,\\
&\text{(Lemma \ref{comm}).}
\end{array}$$
\end{enumerate}

\end{document}